\documentclass{article}
\usepackage[utf8]{inputenc}
\usepackage{epsfig,ulem,amssymb}
\usepackage[margin=0.94in]{geometry}
\usepackage{hyperref,diagbox}
\usepackage{bibunits}
\usepackage{graphicx}
\usepackage{epstopdf}
\usepackage{amsmath}
\usepackage{bm}
\usepackage{setspace}
\usepackage{calligra}

\usepackage{mathrsfs}
\usepackage{subcaption,comment}
\usepackage[dvipsnames]{xcolor}
\usepackage{mathtools}
\usepackage{multirow,enumitem}
\usepackage{lipsum}
\usepackage[affil-it]{authblk}
\usepackage{cite}
\usepackage{leftidx,soul}
\usepackage{relsize}
\usepackage{braket}
\usepackage{relsize}
\usepackage{float}

\usepackage[autostyle]{csquotes}

\usepackage[right]{lineno}
\providecommand{\keywords}[1]{\textbf{\textit{Keywords---}} #1}

\title{Towards a Generalized Approach to Nonlocal Elasticity via Fractional-Order Mechanics}
\author{Sansit Patnaik}
\author{Sai Sidhardh}
\author{Fabio Semperlotti}
\affil{School of Mechanical Engineering, Ray W. Herrick Laboratories, Purdue University, West Lafayette, IN 47907}

\begin{document}
\date{}
\maketitle

\begin{abstract}
This study presents a fractional-order continuum mechanics approach that allows combining selected characteristics of nonlocal elasticity, typical of classical integral and gradient formulations, under a single frame-invariant framework. The resulting generalized theory is capable of capturing both stiffening and softening effects and it is not subject to the inconsistencies often observed under selected external loads and boundary conditions. The governing equations of a 1D continuum are derived by continualization of the Lagrangian of a 1D lattice subject to long-range interactions. This approach is particularly well suited to highlight the connection between the fractional-order operators and the microscopic properties of the medium.
The approach is also extended to derive, by means of variational principles, the governing equations of a 3D continuum in strong form. The positive definite potential energy, characteristic of our fractional formulation, always ensures well-posed governing equations. This aspect, combined with the differ-integral nature of fractional-order operators, guarantees both stability and the ability to capture dispersion without requiring additional inertia gradient terms. The proposed formulation is applied to the static and free vibration analyses of either Timoshenko beams or Mindlin plates. Numerical results, obtained by a fractional-order finite element method, show that the fractional-order formulation is able to model both stiffening and softening response in these slender structures. The numerical results provide the foundation to critically analyze the physical significance of the different fractional model parameters as well as their effect on the  response of the structural elements.\\

\noindent\keywords{Fractional Calculus, Nonlocal Elasticity, Strain Gradient Elasticity, Stiffening, Softening}\\
\noindent All correspondence should be addressed to: \textit{spatnai@purdue.edu} or \textit{fsemperl@purdue.edu}\\

\noindent {\large{\textbf{Highlights}}}
\begin{itemize}
    \item {Fractional-order continuum formulation that captures both stiffening and softening effects.}
    \item {Frame-invariant 3D model developed starting from a 1D lattice with long-range interactions.}
    \item {Well-posed nonlocal governing equations derived from a positive definite system.}
    \item {Predicts anomalous attenuation-dispersion characteristics within a causal framework.}
    \item {Static and free vibration response of Timoshenko beams and Mindlin plates analyzed.}
\end{itemize}
\end{abstract}

%\begin{bibunit}[unsrt]
\section{Introduction}
\label{sec: Introduction}

Several experimental studies have demonstrated that size-dependent effects can become prominent in the response of several structures independently of their spatial scale. In the case of micro- and nano-structures, size-dependent effects have been traced back to material heterogeneity, geometric effects such as changes in curvature, and the existence of surface and interface stresses due to nonlocal atomic interactions and Van der Waals forces \cite{pradhan2009small,wang2011mechanisms,reddy2014non}. Micro- and nano-structures such as carbon nanotubes, thin films and monolayer graphene sheets have far-reaching applications in atomic devices, micro/nano-electromechanical devices, as well as sensors and biological implants. In macroscale applications, particularly those involving heterogeneous structures such as functionally graded materials, metallic foams, granular materials, and porous materials, nonlocal effects have been shown to result from material heterogeneity and interactions between different structural layers  \cite{bazant1976instability,bavzant2000size,hollkamp2019analysis,patnaik2019generalized}. Additionally, specific geometric configurations can also lead to size-dependent effects \cite{buonocore2018tomographic,nair2019nonlocal,buonocore2019occurrence}. In all these macroscopic structures, nonlocal governing equations arise following a homogenization process \cite{bavzant2000size,nair2019nonlocal,buonocore2019occurrence,hollkamp2019analysis}. Based on the examples above, it appears that the ability to accurately model size-dependent effects has profound implications for many engineering applications.

From a general perspective, it is the coexistence of different spatial scales in the above mentioned classes of structural problems that renders the response nonlocal \cite{mindlin1968first,eringen1972linear}. The inability of the classical (i.e. local) continuum theory to capture scale effects prevented its use in these applications and fostered the development of the so-called nonlocal continuum theories. From a general standpoint, the mathematical description of nonlocal continuum theories relies on the introduction of additional contributions in terms of either gradients or integrals of strain (or stress) fields in the constitutive equations. This approach leads to the so-called “weak” gradient methods or “strong” integral methods, respectively. Gradient elasticity theories \cite{mindlin1968first,peerlings2001critical,aifantis2003update,askes2011gradient} account for the nonlocal behavior by introducing strain or stress gradient dependent terms in the stress-strain constitutive law. Integral methods \cite{eringen1972linear,polizzotto2001nonlocal,romano2017stress} capture nonlocal effects by re-defining the constitutive law in the form of a convolution integral of either the strain or the stress field over the horizon of nonlocality. These approaches are further classified as strain-driven or stress-driven \cite{polizzotto2001nonlocal,askes2011gradient,romano2017stress}, depending on whether the nonlocal contributions are modeled using the strain or the stress fields.

Although these different approaches to nonlocal elasticity have been able to address a multitude of aspects typical of the response of size-dependent nonlocal structures, some important challenges still remain open. 
From a high level perspective, gradient theories provide a satisfactory description of the effects of the material microstructure but can introduce significant difficulties connected with the overall stability of the model. As discussed in \cite{askes2011gradient}, while the use of unstable strain-gradients is critical to capture dispersive wave propagation, they give rise to non-convex potential energies leading to the loss of uniqueness in static boundary value problems (BVPs). This issue is often circumvented by using a combination of stable strain-gradients and acceleration gradients \cite{georgiadis2000torsional,metrikine2002one,askes2011gradient}, whose stability comes at the cost of additional terms in both the governing equation and the boundary conditions. From this perspective, integral methods are better suited to deal with boundary conditions and do not lead to any sign paradox, which is peculiar of the gradient methods. However, the corresponding potential energy is not guaranteed to be positive definite and leads to inconsistent predictions for certain loading and boundary conditions \cite{challamel2014nonconservativeness,romano2017stress,romano2017constitutive}. 

From a perspective of practical application, another key limitation of classical nonlocal formulations consists in the fact that, based on the underlying formulation, they can capture only softening or stiffening response but not both simultaneously. Experimental investigations have shown that the size-dependent effects can lead to both stiffening as well as softening of the nonlocal structure depending on the loading and external conditions, such as temperature, loading rate, and boundary conditions \cite{eringen1972linear,bazant1976instability,bavzant2000size,pradhan2009small,wang2011mechanisms,askes2011gradient,karim2013comparison,reddy2014non,guha2015review,fuschi2019size,kloda2020hardening}. To this regard, while classical strain-driven integral formulations \cite{polizzotto2001nonlocal} are suitable for modeling softening effects, stress-driven integral formulations \cite{romano2017stress} and gradient formulations \cite{mindlin1968first} are suitable to capture only stiffening effects. Thus it appears that both the classical integral and gradient formulations are not suitable to capture both stiffening and softening responses. Efforts to achieve an equivalence between the strain-driven integral and gradient formulations, by using special exponential kernels, have been shown to lead to mathematically ill-posed formulations resulting in inaccurate (often called "paradoxical") predictions \cite{romano2017stress,romano2017constitutive}. Further, as stated in \cite{askes2011gradient}, an unresolved issue in strain-gradient formulations pertains to the treatment of materials that exhibit strain-softening. Hence, a comprehensive formulation capable of capturing both stiffening and softening response is still lacking.

In recent years, fractional calculus has emerged as a powerful mathematical tool to model a variety of nonlocal and multiscale phenomena. Fractional derivatives, which are a differ-integral class of operators, are intrinsically multiscale and provide a natural way to account for nonlocal effects. Given the multiscale nature of fractional operators, fractional calculus has found several applications in nonlocal elasticity \cite{atanackovic2009generalized,di2013mechanically,carpinteri2014nonlocal,sumelka2016fractional,sumelka2016geometrical,alotta2017finite,hollkamp2019analysis,patnaik2019generalized,sidhardh2020geometrically,patnaik2020geometrically}. 
Recent studies have shown that a nonlocal continuum approach based on fractional-order kinematic relations provides an effective way to model softening response in nonlocal structures \cite{patnaik2019FEM,patnaik2020geometrically}. These fractional-order nonlocal continuum models result in frame-invariant, thermodynamically consistent and positive definite systems with well-posed governing equations \cite{patnaik2019FEM,patnaik2020plates,sidhardh2020thermoelastic}.

In this study, we show that the differ-integral nature of fractional operators allows them to combine the strengths of both gradient and integral based methods while at the same time addressing a few important shortcomings of both the integer-order formulations. More specifically, we extend the fractional-order continuum formulation developed in \cite{patnaik2019FEM,patnaik2020plates} to develop a comprehensive fractional-order model that captures both softening and stiffening response of nonlocal structures. The overall goal of this study is three fold.

First, we derive the fractional-order governing equations for a 1D nonlocal continuum by continualization of the Lagrangian of a 1D lattice exhibiting long-range interactions with a power-law decay. We will show that fractional-order derivatives of the displacement field (i.e. the nonlocal strain) and fractional-order derivatives of the strain field (i.e. the strain-gradient) are obtained in the potential energy of the 1D structure following continualization of the lattice potential energy. Further, we will demonstrate that the fractional-order formulation is well-posed, frame-invariant, causal, and able to capture anomalous attenuation-dispersion characteristics without the need to resort to acceleration gradient terms, as required in classical strain-gradient formulations. In other terms, in the fractional-order formulation, well-posed governing equations result from a positive definite potential energy while the ability to capture dispersive behavior follows from the differ-integral nature of the fractional operator.
More specifically, the attenuation and dispersion in a solid following the fractional-order formulation are shown to exhibit a power-law dependency on the wave-number/frequency. Remarkably, such anomalous dispersion characteristics have been experimentally observed in different classes of materials including lossy media, fractal and porous materials \cite{szabo1994time,fellah2004verification}, and animal tissues \cite{szabo1994time}. Anomalous attenuation has also been observed in several (non-lossy) scattering media, particularly those characterized by fractal, periodic or random structures \cite{ben2000diffusion,buonocore2018tomographic,buonocore2019occurrence,patnaik2019generalized}. 
Table.~(\ref{tab: comparison}) provides a comparative summary of the classical as well as the fractional-order approaches to nonlocal elasticity, and highlights some of the most distinctive features of the methods.

A second important contribution of this study consists in extending the 1D formulation to a fully 3D formulation. The governing equations in strong form will be derived by using variational principles. In both the 1D and the 3D formulations, we will demonstrate the positive definite nature of the system's potential energy. Additionally, we will discuss the frame-invariance of the formulation and the complete nature of the nonlocal kernel for bounded 3D domains. 

A third key contribution of this work consists in the application of the fractional-order formulation to the analysis of the static and free vibration response of Timoshenko beams and Mindlin plates. The selection of these specific formulations was due to the fact that both the Euler-Bernoulli beam and the Kirchhoff plate formulations can be recovered as special cases; hence, making our study more general and complete. By extending the fractional-order finite element method \cite{patnaik2019FEM,patnaik2020plates} to include the additional gradient terms, we will demonstrate that the fractional-order formulation allows modeling both stiffening and softening effects. We will also critically analyze how the overall structural behavior is affected by the different parameters introduced by the fractional model.

\begin{table}[b]
\centering
\caption{Summary of the fundamental approaches to nonlocal elasticity and comparison of their properties with those of the fractional-order continuum theory. In the table, S.G. denotes strain gradients and I.G. denotes inertia gradients.}
\label{tab: comparison}
\begin{tabular}{|l|c|c|c|c|c|c|}
\hline
\hline
%\multirow{2}{*}{\textbf{Approach type}} &
\multirow{2}{*}{\textbf{\backslashbox{Features}{Approach\\ type}}} &
\multicolumn{2}{c|}{\textbf{Integral}} & \multicolumn{3}{c|}{\textbf{Gradient}} & \multirow{2}{*}{\textbf{Fractional}} \\ \cline{2-6}
 &
  \begin{tabular}[c]{@{}c@{}}Strain\\ driven\end{tabular} &
  \begin{tabular}[c]{@{}c@{}}Stress\\ driven\end{tabular} &
  Stable S.G. &
  Unstable S.G. &
  \begin{tabular}[c]{@{}c@{}}Stable S.G.\\ and I.G.\end{tabular} &
   \\ \hline\hline
\textit{Nature of response}             & Soft          & Stiff         & Stiff    & Stiff    & Stiff   & Soft and Stiff              \\ \hline
\textit{Positive definite system}       & No            & Yes           & Yes      & No       & Yes     & Yes                         \\ \hline
\textit{Capture dispersion}            & Yes           & Yes           & No       & Yes      & Yes     & Yes                         \\ \hline\hline
\end{tabular}
\end{table}

The remainder of the paper is structured as follows: first, we motivate the use of fractional calculus for the analysis of nonlocal structures by considering a 1D lattice with long-range interactions and its corresponding 1D continuum formulation. Next, we extend the 1D continuum to a fully 3D continuum and derive the governing equations in strong form using variational principles. Finally, we use the fractional-order formulation to analyze the effect of the fractional-order nonlocality on the static and free vibration response of beams and plates under different types of loading conditions.

%%%%%%%%%%%%%%%%%%%%%%%%%%%%%%%%%%%%%%%%%%%%%%%%%%%%%%%%%%%%%%%%%%%%%%%%%%%%%%%%%%%%%%%%%%%%%%%%%%%%%%
\section{Fractional-order mechanics: from lattice to 1D continuum}
\label{sec: lattice}
A well established route to develop formulations capable of capturing nonlocal effects in solids is to enforce the continuum limit on a lattice system whose particles are subject to long-range interactions. Several previous works have shown that the continuum limit of lattice structures with one-neighbour and two-neighbour interactions and constant interaction strength lead to the classical first and second integer-order strain-gradient theories of Mindlin, respectively \cite{metrikine2002one,polyzos2012derivation}. An immediate extension of these models follows from considering the response of a lattice with even larger number (i.e. $> 2$) of long-range interactions. Assuming pair-wise constant interaction strengths between different masses across the lattice, it can be easily shown that higher integer-order strain-gradient theories stem from these models. However, these integer-order strain-gradient models would invariably predict a stiffening response of the overall structure. Recall that both softening and stiffening responses have been experimentally observed in the response of solids sensitive to scale effects. In this study, we will show that fractional-order operators can offer a route to develop continuum models capable of predicting both softening and stiffening response in a single formulation. To obtain a physically consistent fractional-order continuum model, we start from a 1D lattice system in which particles are subject to long-range interactions whose pair-wise constant strength decreases with distance in a power-law fashion. While, in the past, other authors have modeled lattices with long-range cohesive forces using fractional calculus \cite{di2013mechanically,carpinteri2014nonlocal}, in this study we extend the formulation by considering also the strain-gradient effects that arise due to microstructural considerations. 

%%%%%%%%%%%%%%%%%%%%%%%%%%%%%%%%%%%%%%%%%%%%%%%%%%%%%%%%%%%%%%%%%%%%%%%%%%%%%%%%%%%%%%%%%%%%%%%%%%%%%%
\subsection{Lattice model and continualization procedure}
\label{ssec: 1D_model}

\begin{figure}[h!]
    \centering
    \includegraphics[width=0.7\textwidth]{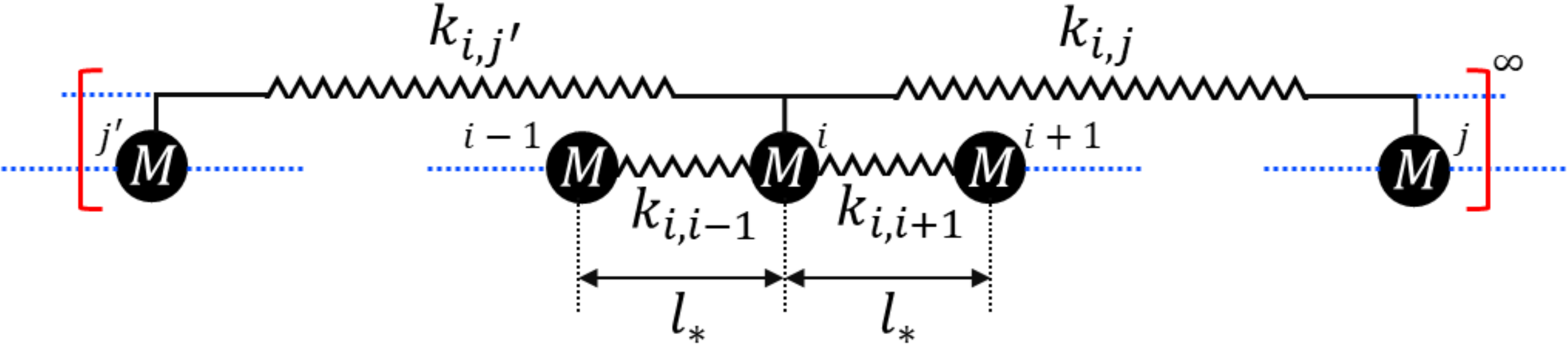}
    \caption{Schematic of the infinite lattice consisting of identical masses denoted as $M$. The masses occur periodically in space separated by a distance of $l_*$. The schematic also illustrates the classical nearest-neighbour interactions as well the long-range interactions between the masses within the infinite lattice.}
    \label{fig: lattice}
\end{figure}

Consider an infinite 1D lattice consisting of identical particles of mass $M$ as shown in Fig.~(\ref{fig: lattice}). The particles are periodically distributed in the $\hat{x}$ direction with spatial period $l_*$ and exhibit only longitudinal motion. The location and displacement of the $n^{th}$ particle (where $n\in \mathbb{Z}$) at the time $t$ are denoted as $x_n(t)$ and $u_n(t)$, respectively. The strength of interaction between particles is modeled via lumped springs having stiffness $k_{i,j}$, where $i^{th}$ and $j^{th}$ are the two interacting particles and $i\neq j$. Note that, in this notation, the comma in the subscript of the spring stiffness does not indicate differentiation. In the following derivation, the dependence of $u_i$ on time $t$ will be implied. Using the above configuration of the lattice and assuming that all the springs are unstressed at the initial time $t=0$, the potential energy stored in the $i^{th}$ cell of the lattice is obtained as:
\begin{equation}
    \label{eq: lattice_PE}
    \mathcal{U}_i = \sum_{j=-\infty}^\infty \frac{1}{2}{k_{i,j}|u_i - u_j|^2}
\end{equation}
where $\mathcal{U}_i$ denotes the potential energy of the $i^{th}$ cell. 
By assuming small displacement gradients (${O}(\varepsilon))$, Taylor's expansion at the point $x_i$ gives:
\begin{equation}
    \label{eq: taylor_series}
    u_i - u_j = (x_i - x_j)\delta_{x_j}^1 u_j + \frac{1}{2}(x_i - x_j)^2 \delta_{x_j}^2 u_j + \textit{h.o.t}
\end{equation}
where $\delta^\square_{x_j}~(\square\in\{1,2\})$ denote the discretized integer-order derivatives at $x_j$.

It is well known that the strength of long-range cohesive forces decays as a function of the inter-atomic distance. Recall that, at continuum level and in integral formulations, this effect is typically accounted for by using convolution terms in the stress-strain constitutive relationships. These convolution kernels have often been chosen to be spatially-decaying exponential functions \cite{eringen1972linear,polizzotto2001nonlocal}. In the lattice model, the stiffness of the springs used to model the interaction between distant particles play a role analogous to the convolution kernels used in classical integral nonlocal elasticity. Thus, in principle, the stiffness of the springs emanating from a given particle towards distant particles can be modeled using spatially decaying exponential functions. In this study, we choose to model the stiffness spatial decay according to power-law functions as follows:
\begin{equation}
    \label{eq: power_law_interaction}
    k_{i,j} =  k_0 \left[\frac{c_1}{|x_{ij}|^{\alpha^2_1} } + \frac{c_2}{|x_{ij}|^{\alpha^2_2} } \right]
\end{equation}
where $|x_{ij}| = |x_{i} - x_{j}|$ indicates the distance between the $i^{th}$ and the $j^{th}$ particles. The parameters $\alpha_1$ and $\alpha_2$ are such that $\alpha_1 \in (0,1)$, $\alpha_2 \in (1,2)$, and $\alpha_2 - \alpha_1 \in (0,1)$. $c_1$ and $c_2$ will be chosen as a function of the parameters $\alpha_1$ and $\alpha_2$, respectively, such that they ensure dimensional consistency and frame-invariance of the formulation. Further, the constant $k_0$ has the dimensions of classical stiffness ($[MT^{-2}]$) and its physical significance will be discussed while deriving the continuum limit of the lattice. Note that the only parameters introduced at this level include $\alpha_1$, $\alpha_2$, and $k_0$. For a given physical lattice with a known spatially decaying stiffness function, these parameters could be obtained by applying standard regression techniques.
Substituting the expression of the stiffness in the infinite series in Eq.~(\ref{eq: lattice_PE}) along with Eq.~(\ref{eq: taylor_series}) and retaining terms up to $O(\varepsilon^2)$, we obtain the potential energy of the $i^{th}$ cell as:
\begin{equation}
    \label{eq: lattice_PE_expansion_step1}
    \mathcal{U}_i = \frac{k_0}{2} \left[ \left(\sum_{j=-\infty}^\infty \frac{ \sqrt{c_1} (x_i - x_j) \delta_{x_j}^1 u_j }{|x_i-x_j|^{\alpha_1}} \right)^2 + \frac{1}{4}\left( \sum_{j=-\infty}^\infty \frac{ \sqrt{c_2} (x_i - x_j) \delta_{x_j}^2 u_j }{|x_i-x_j|^{\alpha_2 - 1}} \right)^2 \right]
\end{equation}

By assuming a small $l_*$ and adopting a continualization process similar to \cite{metrikine2002one,polyzos2012derivation}, the discrete variables indicating the position and the displacement of the particles, can be replaced by the corresponding continuum variables ($x_i \rightarrow x$, $x_j \rightarrow s$, $u_j \rightarrow u(s)$). The constant $k_0$ in the continuum limit can be defined as:
\begin{equation}
    \label{eq: definition_of_k}
    k_0 = \frac{EA}{l_*}
\end{equation}
where $E$ and $A$ denote the Young's modulus and cross-sectional area of the equivalent 1D continuum, respectively. It follows that the constant $k_0$ can be interpreted as the equivalent spring constant representing the strength of the nearest-neighbor interaction forces of a lattice that simulate the microstructure of a local solid (that is not affected by scale effects). Further, we define the constants $c_1$ and $c_2$ in Eq.~(\ref{eq: power_law_interaction}) as:
\begin{equation}
    \label{eq: definition_of_c}
    c_1 = \frac{l_*^2}{4\Gamma(1-\alpha_1)} ~~~~ c_2 = \frac{l_*^4}{4\Gamma(2-\alpha_2)} 
\end{equation}
where $\Gamma(\cdot)$ denotes the Gamma function. %In order to ensure dimensional consistency of the formulation, $c_1$ and $c_2$ are multiplied with the dimensional constants $1\text{m}^{\alpha_1}$ and $1\text{m}^{\alpha_2}$, respectively. 
Under the above assumptions, the continuum limit of the discrete sum in Eq.~(\ref{eq: lattice_PE_expansion_step1}) is obtained to be the following integral representation \cite{ortigueira2006fractional}:
\begin{equation}
    \label{eq: lattice_PE_expansion_step2}
    \mathcal{U}(x) = \frac{EA}{2l_*} \left[ l_*^2 \left[ \frac{1}{2\Gamma(1 - \alpha_1)} \int_{-\infty}^\infty \frac{D_{s}^1 u(s) }{|x-s|^{\alpha_1}} \mathrm{d}s \right]^2 + \frac{l^4_*}{4}\left[ \frac{1}{2\Gamma(2 - \alpha_2)} \int_{-\infty}^\infty \frac{ D_{s}^2 u(s) }{|x-s|^{\alpha_2-1}} \mathrm{d}s
    \right]^2 \right]
\end{equation}
where $D^m_s(\cdot)$ denotes the $m^{th}$ integer-order derivative with respect to the spatial dummy variable $s$ used in the convolution integral.

The convolution integrals in Eq.~(\ref{eq: lattice_PE_expansion_step2}) match with the definition of fractional-order Caputo derivatives with intervals on the real axis, that is $x\in(-\infty,\infty)$ \cite{kilbas2006theory}:
\begin{subequations}
\label{eq: fc_definition}
\begin{equation}
    \label{eq: left}
    {}^{~~~C}_{-\infty} D^{\alpha_m}_{x} u = \frac{1}{\Gamma(m - \alpha_m)} \int_{-\infty}^x \frac{D_{s}^m u(s) }{(x-s)^{\alpha_m-m+1}} \mathrm{d}s
\end{equation}
\begin{equation}
    \label{eq: right}
    {}^{C}_{x} D^{\alpha_m}_{\infty} u = \frac{(-1)^m}{\Gamma(m - \alpha_m)} \int_{x}^\infty \frac{D_{s}^m u(s) }{(s-x)^{\alpha_m-m+1}} \mathrm{d}s
\end{equation}
\end{subequations}
where ${}^{~~~C}_{-\infty} D^{\alpha_m}_{x} (\cdot)$ denotes the left-handed Caputo derivative to the order $\alpha_m$ and lower terminal at $-\infty$, and ${}^{C}_{x} D^{\alpha_m}_{\infty} (\cdot)$ denotes the right-handed Caputo derivative to the order $\alpha_m$ and upper terminal at $\infty$.

Using the above definitions of the left- and right-handed Caputo derivatives, the potential energy density at a point $x$ can be expressed as:
\begin{equation}
    \label{eq: lattice_PE_expansion_step3}
    \Pi(x) = \frac{\mathcal{U}(x)}{Al_*} = \frac{E}{2} \Bigg[ \bigg[ \underbrace{\frac{1}{2} \left( {}^{~~~C}_{-\infty} D^{\alpha_1}_{x} u - {}^{C}_{x} D^{\alpha_1}_{\infty} u \right)}_{\text{Riesz-Caputo derivative}} \bigg]^2 + \frac{l^2_*}{4} \bigg[ \underbrace{\frac{1}{2} \left( {}^{~~~C}_{-\infty} D^{\alpha_2}_{x} u + {}^{C}_{x} D^{\alpha_2}_{\infty} u \right)}_{\text{Riesz-Caputo derivative}} \bigg]^2 \Bigg]
\end{equation}
Recall that, from Eq.~(\ref{eq: power_law_interaction}), $\alpha_1 \in (0,1)$ and $\alpha_2 \in (1,2)$. The above linear combinations of the left- and right-handed Caputo derivatives are typically referred to as the Riesz-Caputo (RC) derivatives. The total potential energy of the structure can now be expressed as:
\begin{equation}
    \label{eq: lattice_PE_expansion_step4}
    \mathcal{U} = \int_{-\infty}^{\infty} \Pi(x)A\mathrm{d}x = \frac{1}{2} \int_{-\infty}^\infty E A \left[ \left( \overline{D}^{\alpha_1}_{x} u \right)^2 + \frac{l^2_*}{4} \left( \overline{D}^{\alpha_2}_x u \right)^2 \right] \mathrm{d}x
\end{equation}
where $\overline{D}^{\alpha_m}_x (\cdot)$ denotes RC derivatives. The over bar $\overline{\square}$ is used to indicate that the RC derivative in Eq.~(\ref{eq: lattice_PE_expansion_step4}) is defined on the real axis, so to differentiate the notation from the RC derivatives defined over bounded domains in \S\ref{sec: 3D continuum}. We merely note that the RC derivative used in the above equation is different from the concept of Riesz derivative defined using sets of Fourier and inverse Fourier transforms \cite{kilbas2006theory}. 

As evident from Eq.~(\ref{eq: lattice_PE_expansion_step4}), the strain in the continuum limit of the infinite lattice structure subject to power-law decaying long-range interactions can be modeled using the RC derivative of the displacement field to the order $\alpha_1 \in (0,1)$. The second term within the integral in Eq.~(\ref{eq: lattice_PE_expansion_step4}) can be interpreted as the fractional-order gradient of the strain field. This is evident by considering the following composition: $\overline{D}^{\alpha_2}_x u = \overline{D}^{\alpha_2-\alpha_1}_x (\overline{D}^{\alpha_1}_x u)$. It follows that we could define a new order $\overline{\alpha}_2=\alpha_2-\alpha_1$. Recall that we have assumed $\alpha_2-\alpha_1\in(0,1)$ in Eq.~(\ref{eq: power_law_interaction}). In order to avoid the introduction of new symbols, we will drop the overline and denote $\overline{\alpha}_2 \equiv \alpha_2$, with the understanding that $\alpha_2$ now lies in the range $(0,1)$. Thus, the total potential energy can be expressed as:
\begin{equation}
    \label{eq: lattice_PE_expansion_step5}
    \mathcal{U} = \frac{1}{2} \int_{-\infty}^\infty E A \Bigg[ \big(\underbrace{\overline{D}^{\alpha_1}_x u }_{\substack{\text{Nonlocal} \\ \text{strain}}} \big)^2 + \frac{l^2_*}{4} \bigg[ \underbrace{\overline{D}^{\alpha_2}_x \left( \overline{D}^{\alpha_1}_x u \right)}_{\substack{\text{Nonlocal gradient} \\ \text{of nonlocal strain}}} \bigg]^2 \Bigg] \mathrm{d}x
\end{equation}
While the specific range for the fractional-orders mentioned here are obtained from mathematical definitions, we will obtain physical constraints on the range of these fractional-orders in \S\ref{ssec: dispersion_1D}.

Given the differ-integral nature of fractional operators, it appears that the different fractional-order derivatives in Eq.~(\ref{eq: lattice_PE_expansion_step5}) lead to a unification of the classical integral and gradient based nonlocal approaches. In fact, the expression in Eq.~(\ref{eq: lattice_PE_expansion_step5}) presents clear insights and comparisons of the fractional-order formulation with both the classical integral and the first-order strain-gradient formulation:
\begin{itemize}
     \item The RC derivative with order $\alpha_1$ captures softening effects in the solid due to the nonlocal interactions. The order $\alpha_1$ captures the strength of the power-law kernel of the fractional derivative which in turn determines the rate of decay in the strength of the nonlocal interactions with distance. Further, the interval of the fractional derivative (here chosen to be $(-\infty,\infty)$), determines the length of the horizon of nonlocality. In other terms, it indicates the distance beyond which nonlocal interactions are no longer accounted for in the fractional derivative \cite{patnaik2019FEM,patnaik2020plates}.

     \item From Eqs.~(\ref{eq: lattice_PE_expansion_step1},\ref{eq: lattice_PE_expansion_step2}) it is seen that, for the lattice with long-range cohesive interactions, the expression for the potential energy at a point $x$ includes contribution of the microstructural information (that is the strain-gradient) of all points in the nonlocal horizon of $x$. This is in addition to the nonlocal contribution of the strain energy captured by the RC derivative $\overline{D}^{\alpha_1}_x u$. It is immediate to see that the RC derivative of the nonlocal strain with order $\alpha_2$ captures the stiffening effects in the solid. More specifically, analogous to classical strain-gradient formulations, this term would account for the microstructural information within the strain energy potential. Furthermore, the parameter $l_*$ that was initially introduced as the lattice parameter can be interpreted as the microstructural length scale analogously to classical formulations.
\end{itemize}
The above discussions lead to the conclusion that the use of the different fractional-order gradients allows the continuum model to capture simultaneously both long-range cohesive forces (leading to softening effects) as well as strain-gradient terms capturing microstructural properties (leading to stiffening effects). A remarkable outcome of this approach is that, not only it can capture both softening and stiffening effects in a single formulation, but it can account for these effects simultaneously. Note that the first-order strain-gradient theory for the 1D continuum can be obtained from the above formulation by using $\alpha_1=1$ and $\alpha_2=1$. Following the above discussion, we call $\alpha_1$ as \textit{nonlocal-strain order} and $\alpha_2$ as the \textit{strain-gradient order}.

Note that the definition of the spring stiffness in Eq.~(\ref{eq: power_law_interaction}) leads to $k_{i,j}=k_{j,i}$. This ensures that the internal state of the lattice cannot be changed following a translation of all the particles by the same distance. While this is sufficient to ensure frame-invariance of the 1D continuum, the extension to a full 3D model would require the satisfaction of frame-invariance under rotations as well. It is also important to note that the potential energy of the nonlocal 1D solid consists of Caputo derivatives and not other types of fractional derivatives (e.g. Riemann Liouville). Recall that the Caputo derivative of a constant function is zero, as for classical integer order derivatives. This property does not hold true for all definitions of fractional derivatives \cite{kilbas2006theory}. However, in the context of frame invariance, this is a key point that ensures that no strain is accumulated in the 1D solid under translation, that is for a constant $u(x)$.

The kinetic energy of the 1D solid can be evaluated similar to classical integer-order formulations. Note that the introduction of nonlocality through the long-range spring connections has no effect on the expression for kinetic energy. It follows that, in the continuum limit, the kinetic energy of the above described 1D solid is given as \cite{metrikine2002one,polyzos2012derivation}:
\begin{equation}
    \label{eq: lattice_KE}
    T = \frac{1}{2} \int_{-\infty}^{\infty} \left[ \rho A (D^1_t u)^2 + \rho^\prime A \frac{l^2_*}{3} \left[ D^1_x \left( D^1_t u \right) \right]^2 \right]  \mathrm{d}x
\end{equation}
where $D^1_t(\cdot)$ denotes the first integer-order derivative with respect to time, $\rho$ is the density, and $\rho^\prime$ is the microdensity of the solid that has the same interpretation as in classical integer-order strain-gradient models. A possible extension of the fractional-order continuum theory developed above involves the use of time fractional derivatives within the kinetic energy as:
\begin{equation}
    \label{eq: fractional_lattice_KE}
    T = \frac{1}{2} \int_{-\infty}^{\infty} \left[ \rho A ({}^C_0 D^\kappa_t u)^2 + \rho^\prime A \frac{l^2_*}{3} \left[ D^1_x \left( {}^C_0 D^\kappa_t u \right) \right]^2 \right]  \mathrm{d}x
\end{equation}
where ${}^C_0 D^\kappa_t u$ is a left-handed Caputo derivative with order $\kappa\in(0,1)$ and defined on the interval $(0,t)$. This will allow the fractional-order model to capture memory effects and non-conservative dissipation mechanisms, such as those encountered in viscoelastic materials. Such a formulation can be found in \cite{ansari2017studying} where the nonlinear response of viscoelastic nanobeams have been captured by using time fractional derivatives. However, unlike our study, size-dependent effects in \cite{ansari2017studying} were modeled using the classical first-order strain-gradient formulation. Since memory effects and dissipation have already been addressed in the literature, in this study we focus on the modeling of nonlocal effects in non-dissipative solids using space fractional derivatives.

%%%%%%%%%%%%%%%%%%%%%%%%%%%%%%%%%%%%%%%%%%%%%%%%%%%%%%%%%%%%%%%%%%%%%%%%%%%%%%%%%%%%%%%%%%%%%%%%%%%%%%
\subsection{Governing equations for the 1D continuum}
\label{ssec: 1D_governing_equations}
We derive the dynamic governing equations of the 1D structure by using Hamilton's variational principle:
\begin{equation}
    \label{eq: Hamilton_eq}
    \int_{t_0}^{t_1} \delta(\mathcal{U} - T)\mathrm{d}t = \int_{t_0}^{t_1} \delta\left[ \frac{A}{2} \int_{-\infty}^\infty \Big[ E(\overline{D}^{\alpha_1}_x u )^2 + \frac{l^2_*}{4} E \left[ {\overline{D}^{\alpha_2}_x \left( \overline{D}^{\alpha_1}_x u \right)} \right]^2 -\rho (D^1_t u)^2 - \rho^\prime \frac{l^2_*}{3} \left[ D^1_x \left( D^1_t u \right) \right]^2 \Big] \mathrm{d}x \right] \mathrm{d}t 
\end{equation}
Performing variational simplifications, the governing equation is obtained as:
\begin{equation}
    \label{eq: 1D_gov_equations}
    E \left[ \overline{D}^{2\alpha_1}_x u - \frac{l^2_*}{4} \overline{D}^{2(\alpha_1+\alpha_2)}_x u \right] = \rho D^2_t u - \rho^\prime \frac{l^2_*}{3} D^2_x \left( D^2_t u \right)
\end{equation}
A detailed derivation of the above equations is provided for a 3D bounded continuum in the Supplementary Information (SI).

Recall that capturing dispersive wave propagation is one of the main motivation promoting the development of gradient elasticity in classical elastodynamics. As discussed in detail in \cite{askes2011gradient}, the use of "unstable" (integer-order) strain-gradients is critical in capturing wave dispersion, however, in the static sense, "unstable" strain-gradients result in non-convex potential energies leading to the loss of uniqueness in static boundary value problems (BVPs).
In the classical analogue of Eq.~(\ref{eq: 1D_gov_equations}), a positive (negative) sign of the strain-gradient term corresponds to an unstable (stable) strain-gradient. While the combined used of these gradient terms is generally avoided because one of the two terms will always tend to predominate, this issue is circumvented by using a combination of stable (integer-order) strain-gradients and acceleration gradients (see, \cite{georgiadis2000torsional,metrikine2002one}) which allows for dispersive wave propagation while ensuring a well-posed BVP. A detailed discussion on this aspect can be found in \cite{askes2011gradient}, where a combination of different stain and acceleration gradients \footnote{Different researchers have used different terminology (acceleration-gradient or velocity-gradient) to refer to the term $D^2_x (D^2_t u)$. We be believe that both the terminology are appropriate since the term appears as an acceleration gradient in the strong form and translates to a velocity gradient in weak form. In this study, following \cite{askes2011gradient}, we refer to it as the acceleration gradient.} is studied to arrive at theories which are well suited for both static and dynamic applications. 
To this regard, we highlight that the fractional-order strain-gradient formulation provides a natural way of dealing with this issue without the need of additional stabilising acceleration gradients. Note that the potential energy given in Eq.~(\ref{eq: lattice_PE_expansion_step5}), resulting from the fractional-order formulation, is quadratic in nature and hence fully convex. Additionally, it is established in \cite{patnaik2019FEM} that the fractional-order operators are self-adjoint and the resulting formulation leads to well posed BVPs. 
Further, the specific form of the spring strength given in Eq.~(\ref{eq: power_law_interaction}) indicates that the stiffness of the structure exhibits dependence on wavelength and hence, the fractional-order formulation, obtained via continualization of the Lagrangian of the 1D lattice, is well suited to capture anomalous dispersion characteristics (\S\ref{ssec: dispersion_1D}). Further, we will establish in the following \S\ref{ssec: dispersion_1D} that the fractional-order formulation is causal and stable.\\

%%%%%%%%%%%%%%%%%%%%%%%%%%%%%%%%%%%%%%%%%%%%%%%%%%%%%%%%%%%%%%%%%%%%%%%%%%%%%%%%%%%%%%%%%%%%%%%%%%%%%%
\subsection{Dispersion analysis of the 1D continuum}
\label{ssec: dispersion_1D}
To obtain the dispersion relation, we substitute in the fractional-order elastodynamic equation given in Eq.~(\ref{eq: 1D_gov_equations}) the following ansatz:
\begin{equation}
    \label{eq: 1D_solution_assumption}
    u(x,t) = u_0 e^{i(kx - \omega t)}
\end{equation}
where $u_0$ is the amplitude of the longitudinal wave, $k$ denotes the wave-number, $\omega$ denotes the angular frequency of free longitudinal vibrations, and $i=\sqrt{-1}$. For the RC derivatives on the real line used in Eq.~(\ref{eq: fractional_lattice_KE}) \cite{kilbas2006theory}:
\begin{equation}
    \label{eq: fd_of_exp}
    D^\alpha_x (e^{kx}) = k^\alpha e^{kx}
\end{equation}
Using the above RC derivative of the exponential, we obtain the complete form of the dispersion relations for longitudinal waves in the 1D solid as:
\begin{equation}
    \label{eq: 1D_dispersion}
    \frac{\omega}{k} = \sqrt{\frac{E}{\rho}} {\left[ -i^{2\alpha_1} k^{2(\alpha_1 - 1)} + i^{2(\alpha_1 + \alpha_2)} k^{2(\alpha_1 + \alpha_2-1)}\frac{l^2_*}{4} \right]^\frac{1}{2}} \left[ 1+\frac{\rho^\prime l^2_*}{3\rho} k^2 \right]^{-1}
\end{equation}
Using Euler's formula, the above equation can be recast in the following manner:
\begin{equation}
    \label{eq: 1D_dispersion_final}
    \begin{split}
    \frac{\omega}{k} = \mathcal{Z} = \Bigg[ \bigg( &-\cos(\alpha_1\pi) k^{2(\alpha_1 - 1)} +  \cos({2(\alpha_1 + \alpha_2)}\pi) k^{2(\alpha_1 + \alpha_2-1)}\frac{l^2_*}{4} \bigg)  +  \\ &  i \bigg(\underbrace{-\sin(\alpha_1\pi) k^{2(\alpha_1 - 1)} + \sin({2(\alpha_1 + \alpha_2)}\pi) k^{2(\alpha_1 + \alpha_2-1)}\frac{l^2_*}{4}}_{\overline{b}} \bigg) \Bigg]^\frac{1}{2} {\underbrace{\left[ 1+\frac{\rho^\prime l^2_*}{3\rho} k^2 \right]}_{\phi}}^{-1}
    \end{split}
\end{equation}

Expressing $\omega=\mathcal{Z}k$, stable and causal solutions are recovered when $\Re(\mathcal{Z})>0$ and $\Im(\mathcal{Z})<0$. Note from Eq.~(\ref{eq: 1D_solution_assumption}) that $\Re(\mathcal{Z})>0$ would lead to forward propagating solutions ensuring causality, while $\Im(\mathcal{Z})<0$ leads to attenuation hence ensuring stability. Thus, it appears that the complex number $\mathcal{Z}$ must lie in the fourth quadrant of the Argand plane or, equivalently, $\mathcal{Z}^2$ must lie below the real-axis of the Argand plane. It immediately follows that the quantity $\overline{b} = \Im(\mathcal{Z}^2)$ in Eq.~(\ref{eq: 1D_dispersion}) must be less than or equal to zero for all values of $k$ and $l_*$. The latter condition holds true for all positive values of the wave-number $k$ and microstructural length $l_*$ under the following restrictions for $\alpha_1$ and $\alpha_2$:
\begin{equation}
    \label{eq: restriction_on_order}
    \alpha_1, \alpha_2 \in [0.5,1] 
\end{equation}
Under the above condition, $\sin(\alpha_1 \pi) > 0$ and $\sin ({2(\alpha_1 + \alpha_2)}\pi) <0$, ensuring that $\overline{b} < 0$ for all positive values of $k$ and $l_*$.
It follows that, in this study, we only consider values of the fractional-orders which lie in $[0.5,1]$. Under the above conditions, the $\Re(\mathcal{Z})$ would contribute to anomalous wave-number dependent dispersion in the propagating longitudinal waves while $\Im(\mathcal{Z})$ would lead to attenuation in the propagating waves. 

Note that the term indicated by $\phi$ in Eq.~(\ref{eq: 1D_dispersion_final}) appears from the inclusion of the acceleration gradient term in the governing equations. As discussed in \cite{askes2011gradient,sidhardh2019dispersion}, the inclusion of the acceleration gradient term prevents an unbounded growth in the wave speed following an increase in the wave number. We merely note that, given the attenuation in the wave speed, the inclusion of the acceleration gradient term is no longer necessary in the fractional-order formulation. To this regard, note that ignoring the term $\phi$ would cause the dispersion as well as the attenuation in the longitudinal wave speeds to exhibit a power-law dependence on the wave-number. This is a direct consequence of the power-law nature of the strength of the long-range interactions. Remarkably, several studies have highlighted a power-law dependence of the attenuation-dispersion relations on frequency/wave-number in many types of lossy and highly scattering media, including fractal and porous materials, and animal tissues \cite{szabo1994time,ben2000diffusion,fellah2004verification}. It follows that, in this study, we will neglect the acceleration gradients and focus on modeling media with power-law attenuation-dispersion behavior.
Another particularly interesting outcome of the above formulation is that the dispersion and attenuation form a Hilbert pair, ensuring that the dynamic formulation is fully causal \cite{szabo1994time,fellah2004verification,patnaik2019generalized}.

%%%%%%%%%%%%%%%%%%%%%%%%%%%%%%%%%%%%%%%%%%%%%%%%%%%%%%%%%%%%%%%%%%%%%%%%%%%%%%%%%%%%%%%%%%%%%%%%%%%%%%
\section{Extension to 3D continuum}
\label{sec: 3D continuum}
The previous section used a 1D framework to illustrate the remarkable features of the fractional-order formulation. In this section, we extend the formulation to a fully three-dimensional and finite solid. The governing equations for the 3D continuum are derived using Hamilton's variational principle. We highlight here that, although the 3D formulation presented in the following is developed by continualization of the 1D lattice, the same formulation can also be derived from a continuum-mechanics approach by considering different configurations of a nonlocal solid, as illustrated in \cite{patnaik2019generalized}. More specifically, the 3D formulation developed in this study via continualization principles can also be obtained from the fractional-order continuum formulation presented in \cite{patnaik2019generalized} by adding fractional-order strain-gradient terms to the constitutive relations. To this regard, note that the continualization route adopted in this study motivates the need of a fractional-order approach to capture both stiffening and softening effects within a single formulation.

%%%%%%%%%%%%%%%%%%%%%%%%%%%%%%%%%%%%%%%%%%%%%%%%%%%%%%%%%%%%%%%%%%%%%%%%%%%%%%%%%%%%%%%%%%%%%%%%%%%%%%
\subsection{Weak formulation}
\label{ssec: weak_formulation_3D}
The potential energy derived for the nonlocal 1D continuum in Eq.~(\ref{eq: lattice_PE_expansion_step5}) is extended to a 3D continuum in the following manner:
\begin{equation}
\label{eq: 3D_PE}
    \mathcal{U} = \frac{1}{2}\int_{\Omega} \left[ {\bm\epsilon} : \bm{C} : {\bm\epsilon} + {\bm\eta}: \bm{G} : {\bm\eta} \right] \mathrm{d}\mathbb{V}
\end{equation}
where $\bm{C}$ denotes the classical fourth-order elasticity tensor and $\bm{G}$ is the sixth-order elasticity tensor. ${\bm\epsilon}$ and ${\bm\eta}$ denote the fractional-order strain and its gradient, respectively. The volume of the 3D continuum is denoted by $\Omega$ and $\mathrm{d}\mathbb{V}$ denotes an infinitesimal volume element. Note that the total potential energy is positive definite for positive definite material elasticity tensors. 

\begin{figure}[h]
	\centering
	\includegraphics[width=\linewidth]{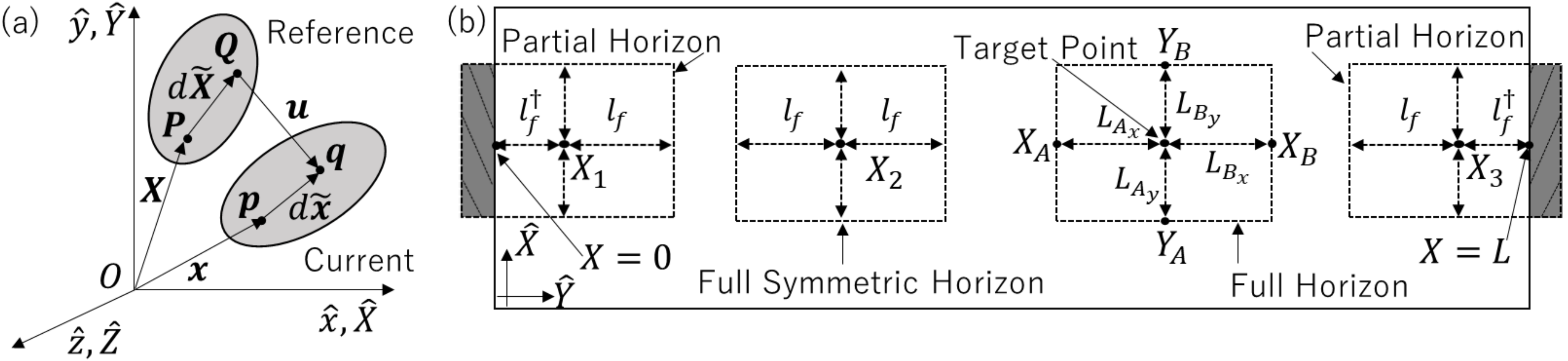}
	\caption{\label{fig: FCM} (a) Schematic indicating the infinitesimal material $\mathrm{d}\tilde{\bm{X}}$ and spatial $\mathrm{d}\tilde{\bm{x}}$ line elements in the nonlocal medium subject to the displacement field $\bm{u}$. (b) Horizon of nonlocality and length scales at three different material points $\bm{X}_1$, $\bm{X}_2$, and $\bm{X}_3$ in a 2D domain. Note that in the $\hat{X}$ direction, $\bm{X}_2$ has a horizon of nonlocality equal to $l_f$ on both the left and the right sides, while the horizon of nonlocality at the points $\bm{X}_1$ and $\bm{X}_3$ are truncated to $l_f^\dagger$ such that $l_f^\dagger<l_f$, on the left and the right sides, respectively. Clearly, the nonlocal model can account for a partial (i.e. asymmetric) horizon condition that occurs at points close to a boundary or interface.}
\end{figure}

The infinitesimal strain in the 3D nonlocal continuum is obtained by extending the 1D nonlocal strain indicated in Eq.~(\ref{eq: lattice_PE_expansion_step5}) as:
\begin{equation}
\label{eq: infinitesimal_fractional_strain}
     {\bm{\epsilon}}=\frac{1}{2} \bigr(\bm\nabla^{\alpha_1} {\bm{U}}_X + \bm\nabla^{\alpha_1} {\bm{U}}_X^{T}\bigl)=\frac{1}{2} \bigr(\bm\nabla^{\alpha_1} {\bm{u}}_x + \bm\nabla^{\alpha_1} {\bm{u}}_x^{T}\bigl)
\end{equation}
where $\bm{U}(\bm{X})=\bm{x}(\bm{X})-\bm{X}$ and $\bm{u}(\bm{x})=\bm{x}-\bm{X}(\bm{x})$ are the displacement fields in the Lagrangian ($\bm{X}$) and Eulerian ($\bm{x}$) coordinates, respectively (see Fig.~(\ref{fig: FCM})a). $\bm{\nabla}^{\alpha_m}(\cdot)~(\alpha_m\in\{\alpha_1,\alpha_2\})$ is the RC fractional gradient operator defined as:
\begin{equation}
    \label{eq: RC_grad_operator}
    \bm{{\nabla}}^{\alpha_m} (\cdot) = {D}_{x}^{\alpha_m} (\cdot) \hat{x} + {D}_{y}^{\alpha_m} (\cdot) \hat{y} + {D}_{z}^{\alpha_m} (\cdot) \hat{z}
\end{equation}
where $\{\hat{x},\hat{y},\hat{z}\}$ are the Cartesian basis vectors. $D^{\alpha_m}_{x_j}(\cdot)$ are the RC fractional derivatives which will be defined in the following.
We emphasize that the above definition for the strain tensor can also be derived rigorously following a continuum mechanics approach, starting from a fractional-order definition of the deformation gradient tensor (see \cite{patnaik2019generalized,patnaik2019FEM,sumelka2016fractional}). Further, the fractional gradient of the nonlocal strain is defined as:
\begin{equation}
	\label{eq: frac_straingrad}
	{\bm\eta}=\bm\nabla^{\alpha_2}{\bm\epsilon}
\end{equation}
It follows that the constitutive relations for the Cauchy stress and the higher-order stress, in terms of the work-conjugates ${\bm\epsilon}$ and $\bm{\eta}$, can be expressed as:
\begin{subequations}
\label{eq: constt}
\begin{equation}
	{\bm\sigma} = \bm{C} : \bm\epsilon
\end{equation}
\begin{equation}
   \bm{\tau} = \bm{G} : \bm{\eta}
\end{equation}
\end{subequations}

While the RC fractional derivatives used for the infinite 1D solid in \S\ref{sec: lattice} were defined on the real axis, these derivatives are modified for bounded domains to ensure frame-invariance everywhere on the domain and a complete kernel when approaching boundaries \cite{patnaik2019generalized,patnaik2019FEM}. Note that completeness of the kernel in nonlocal elasticity is critical to ensure well-posed problems and stable numerical implementations. The space-fractional derivative $D^{\alpha_m}_{\bm{X}} \bm\psi(\bm{X},t)$ $(m\in\{1,2\})$ of the function $\bm\psi(\bm{X},t)$ $(=\bm{U}(\bm{X},t)$ or ${\bm\epsilon(\bm{X},t)})$ in Eqs.~(\ref{eq: infinitesimal_fractional_strain},\ref{eq: frac_straingrad}) is taken according to a RC definition with order $\alpha_m\in(0.5,1)$ defined on the interval $\bm{X} \in (\bm{X}_A,\bm{X}_B) \in \mathbb{R}^3 $. The RC definition for this bounded domain is defined as a linear combination of the left- and right-handed Caputo derivatives in the following manner \cite{patnaik2019generalized}:
\begin{subequations}
	\label{eq: RC_definition}
	\begin{equation}
	D^{\alpha_m}_{\bm{X}} \bm\psi(\bm{X},t) = \frac{1}{2}\Gamma(2-\alpha_m) \big[\bm{L}_{A}^{\alpha_m-1}~ {}^{~~C}_{\bm{X}_{A}}D^{\alpha_m}_{\bm{X}} \bm\psi (\bm{X},t) - \bm{L}_{B}^{\alpha_m-1}~ {}^C_{\bm{X}}D^{\alpha_m}_{\bm{X}_{B}} \bm\psi (\bm{X},t) \big]
	\end{equation}
	\begin{equation}
	D^{\alpha_m}_{X_j} \psi_i(\bm{X},t) = \frac{1}{2}\Gamma(2-\alpha_m) \big[L_{A_j}^{\alpha_m-1}~{}^{~~~C}_{X_{A_j}} D^{\alpha_m}_{X_j} \psi_i(\bm{X},t) - L_{B_j}^{\alpha_m-1}~{}^{~C}_{X_j}D^{\alpha_m}_{X_{B_j}} {\psi_i}(\bm{X},t)\big]
	\end{equation}
\end{subequations}
where, ${}^{~~C}_{\bm{X}_{A}}D^{\alpha_m}_{\bm{X}} \bm\psi(\bm{X},t)$ and ${}^C_{\bm{X}}D^{\alpha_m}_{\bm{X}_{B}} \bm\psi(\bm{X},t)$ are the left- and right-handed Caputo derivatives of $\bm\psi(\bm{X},t)$ respectively. In the indicial expression in Eq.~(\ref{eq: RC_definition}b), $L_{A_j}$ and $L_{B_j}$ are length scales along the $j^{th}$ direction in the reference configuration. 
The index $j$ in Eq.~(\ref{eq: RC_definition}b) is not a repeated index because the length scales are scalar multipliers. In the current configuration, these length scales are denoted as $l_{A_j}$ and $l_{B_j}$. The interval of the fractional derivative $(\bm{X}_A,\bm{X}_B)$ defines the horizon of nonlocality which is schematically shown in Fig.~(\ref{fig: FCM}) for a generic point $\bm{X}\in\mathbb{R}^2$.
This interval defines the set of all points in the solid that influence the elastic response at $\bm{X}$ or, equivalently, the characteristic distance beyond which information of nonlocal interactions is no longer accounted for in the derivative. 
Recall that the use of Caputo derivatives ensured a frame-invariant model for the 1D continuum. As discussed in \cite{patnaik2019generalized}, the terms $\frac{1}{2}\Gamma(2-\alpha_m)$, $L_{A_j}^{\alpha_m-1}$, and $L_{B_j}^{\alpha_m-1}$ ensure the frame invariance of the 3D formulation. Further, it is required that the length scales $\bm{L}_A=\bm{X}-\bm{X}_A$ and $\bm{L}_B=\bm{X}_B-\bm{X}$. Hence, it follows that the length scales, $L_{A_j}$ and $L_{B_j}$ physically denote the dimension of the horizon of nonlocality to the left and right of point $\bm{X}$ along the $j^{th}$ direction. The length scales have been schematically illustrated in Fig.~(\ref{fig: FCM}b). The introduction of the different length scales ($\bm{L}_{A}$ and $\bm{L}_B$) is to enable the formulation to deal with possible asymmetries in the horizon of nonlocality (e.g. resulting from a truncation of the horizon when approaching a boundary or an interface). Note also that the length scale parameters ensure the dimensional consistency of the formulation.

A key aspect in nonlocal integral formulations is the nature of the kernel when approaching the boundaries. To this regard, we highlight that the definition of the RC derivative in Eq.~(\ref{eq: RC_definition}) ensures the completeness of the power-law convolution kernel within the fractional-order derivative. Note that the lower terminal is $\bm{X}_A=\bm{X}-\bm{L}_A$ and the upper terminal is $\bm{X}_B=\bm{X}+\bm{L}_B$. This definition allows the length scales $\bm{L}_A$ and $\bm{L}_B$ to be truncated when the point $\bm{X}$ approaches a boundary (see Fig.~(\ref{fig: FCM}b)). It follows that the terminals of the RC derivative are properly modified hence resulting in a complete kernel over the truncated domain. The completeness of the kernel can also be established by investigating the nature of the fractional-order model at points on the boundary, that is when either $L_{A_j}\rightarrow0$ or $L_{B_j}\rightarrow0$. As established in \cite{patnaik2019generalized,patnaik2019FEM}, for a material point (say $\bm{X}_0$) located on one of the boundaries (identified by the normal in the $j^{th}$ direction), for the limiting case when $L_{A_j}\rightarrow0$, the RC fractional derivative reduces to:
\begin{equation}
\label{sq13}
\lim_{L_{A_j}\to 0} D^{\alpha_m}_{X_j} \psi_i (\bm{X},t)= \frac{1}{2} \left[ \frac{\mathrm{d} \psi_i(\bm{X},t)}{\mathrm{d}X_j} \bigg|_{\bm{X}_0} + (1-\alpha_m) L_{B_{j}}^{\alpha_m-1} \int_{X_{0_j}}^{X_{B_j}} \frac{D^1_{S_j}\psi_{i}(\textbf{S},t)}{(S_j-X_j)^{\alpha_m}} dS_j \right]
\end{equation}
where $S_j$ is a dummy vector variable used to carry out the spatial convolution integral. From Eq.~(\ref{sq13}) it is immediate to observe that while the right-handed Caputo derivative captures nonlocality ahead of the point $X_0$ (in the $j^{th}$ direction), the left-handed derivative is reduced to the classical first-order derivative. This suggests that the truncation of the nonlocal horizon (and the corresponding convolution) at the boundary has been accounted for in a consistent manner. Similar expressions hold when $L_{B_j}=0$ and for the deformed configuration ($l_{A_j}=0$ or $l_{B_j}=0$). 

The above discussions on the frame-invariance of the formulation and on the nature of the kernel close to material boundaries establish both the completeness and consistency of the fractional-order continuum formulation. It remains to obtain the expressions for the kinetic energy of the continuum and the work done by externally applied forces. The work done by external forces is defined analogous to classical formulations of gradient elasticity as:
\begin{equation}
\label{eq: 3D_ExtW}
   \mathcal{V} = \int_{\Omega} ({\overline{\bm{b}}} \cdot \bm{u} ) \mathrm{d} \mathbb{V} +\int_{\partial \Omega} ({\overline{\bm{t}}} \cdot \bm{u} +\bm{\overline{q}} \cdot \bm{\hat{n}} \cdot ( \bm{u} \otimes \bm{\nabla}^{\alpha_1})) \mathrm{d} \mathbb{A} + \oint_{\overline{\Gamma}} (\bm{\overline{r}} \cdot \bm{u}) \mathrm{d}\mathrm{l}
\end{equation}
where $\mathrm{d}\mathrm{A}$ and $\mathrm{d}\mathrm{l}$ indicate area and line elements along the surface $\partial \Omega$ (with normal $\bm{\hat{n}}$) and edge $\overline{\Gamma}$ of the solid, respectively. The bar on $\overline{\Gamma}$ symbol in the above equation, is used to differentiate the same from the $\Gamma(\cdot)$ function and the symbol $\otimes$ denotes the dyadic product. $\bm{\bar{b}}$, $\bm{\bar{t}}$, $\bm{\bar{q}}$, and $\bm{\bar{r}}$ are the prescribed values of body force per unit volume, surface traction per unit area, double stress traction vector and line load along sharp edges of the continuum, respectively. Finally, recalling that the introduction of nonlocality has no effect on the expression of the kinetic energy, we can write:
\begin{equation}
\label{eq: 3D_KE}
    T = \frac{1}{2} \int_{\Omega} \rho (\bm{\dot{u}} \cdot {\bm{\dot{u}}})\mathrm{d}\mathbb{V} 
\end{equation}
where $\rho$ indicates the density of the solid and $\dot{\square}$ indicates the first integer-order derivative with respect to time. By using the Hamilton's principle and the expressions of the potential energy, kinetic energy, and work done by external forces, the weak form of the governing equations for the 3D continuum are expressed as:
\begin{equation}
    \label{eq: extended_hamiltons_principle}
    \int_{t_1}^{t_2}\left(\delta \mathcal{U} - \delta V - \delta T\right)\mathrm{d}t=0
\end{equation}

%%%%%%%%%%%%%%%%%%%%%%%%%%%%%%%%%%%%%%%%%%%%%%%%%%%%%%%%%%%%%%%%%%%%%%%%%%%%%%%%%%%%%%%%%%%%%%%%%%%%%%
\subsection{Strong formulation}
\label{ssec: strong_formulation_3D}
The strong form of the fractional-order governing equations are obtained by applying the fundamental law of variational calculus to Eq.~(\ref{eq: extended_hamiltons_principle}). Analogously to classical integer-order formulations, the procedure to obtain the strong form for 3D domains involves the use of different principles of vector calculus. 
To this regard, note that fractional vector calculus principles have been recently developed and do not hold true for a general bounded geometry \cite{tarasov2008fractional}. This aspect can be attributed to the fact that fractional-order operators (i.e. derivatives or integrals) do not generally commute, except when defined on the real axis \cite{kilbas2006theory,tarasov2008fractional}. However, we will show that the variational statement in Eq.~(\ref{eq: extended_hamiltons_principle}) can be exactly simplified when considering a cuboidal (or, rectangular) geometry. It can also be envisioned that, the strong form derived assuming a cuboidal geometry will also be applicable for geometries wherein the surfaces/edges can be exactly represented or even approximated by using rectangular/line elements. Although the strong form requires the simplified cuboidal geometry, we emphasize that the weak form in Eq.~(\ref{eq: extended_hamiltons_principle}) is applicable to any geometry.

\begin{figure}[h]
    \centering
    \includegraphics[width=0.3\textwidth]{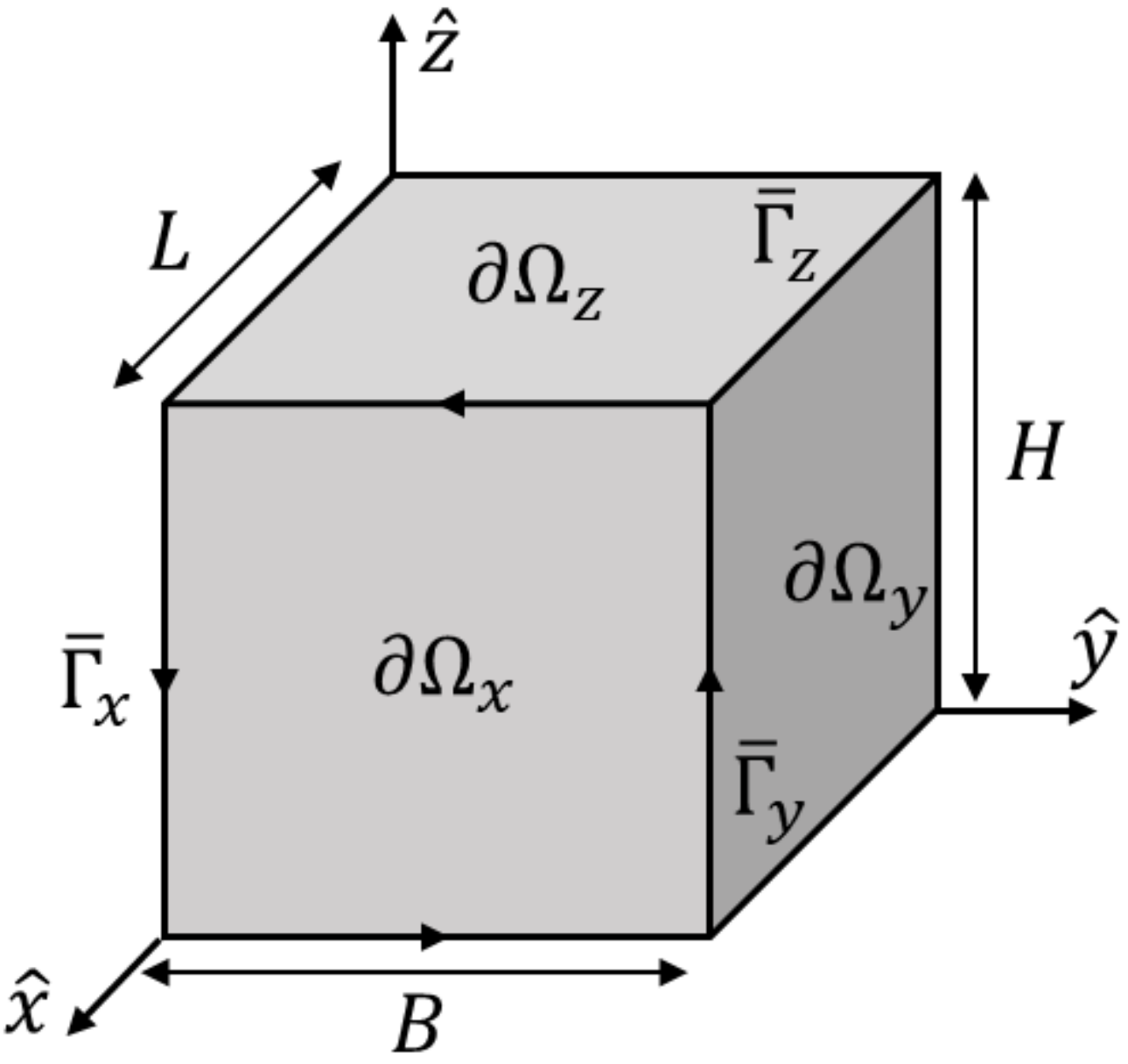}
    \caption{Schematic of the cuboidal domain ($\Omega$) illustrating the different geometrical parameters. The surface of the cuboid is given as $\partial \Omega = \partial \Omega_x \cup \partial \Omega_y \cup \partial \Omega_z$, where $\partial\Omega_{x_k}$ denotes a surface with its normal oriented along the positive or negative $\hat{x}_k$ axis. The edges of the cuboid are denoted by $\overline{\Gamma}=\overline{\Gamma}_x \cup \overline{\Gamma}_y \cup \overline{\Gamma}_z$. $\overline{\Gamma}_{x_k}$ denotes the edges of the surface $\partial \Omega_{x_k}$ oriented in the anti-clockwise sense with respect to the normal to the surface.}
    \label{fig: cube}
\end{figure}

Considering the cuboidal geometry $\Omega$ illustrated in Fig.~(\ref{fig: cube}), the first variation of the potential energy is obtained as:
\begin{equation}
\label{eq: PE_variation}
\begin{split}
\delta \mathcal{U} = -\int_{\Omega} {\bm{\widetilde{\nabla}}}^{\alpha_1} \cdot \big( \bm{\sigma} - {\bm{\widetilde{\nabla}}}^{\alpha_2} \cdot \bm{\tau} \big) \cdot \delta \bm{u} \mathrm{d}\mathbb{V} 
+ \int_{\partial \Omega} \left[ \bm{I}^{1-\alpha_1}_{\bm{\hat{n}}} \cdot \big( \bm{\sigma} - {\bm{\widetilde{\nabla}}}^{\alpha_2} \cdot \bm{\tau} \big) - \left[ \bm{R} : \left( \bm{I}^{1-\alpha_2}_{\bm{\hat{n}}} \cdot \bm{\tau} \right) \otimes {\bm{\widetilde{\nabla}}}^{\alpha_1} \right] \right] \cdot \delta \bm{u} \mathrm{d}\mathbb{A}  + \\
+\int_{\partial \Omega} \left[ \bm{I}^{1-\alpha_2}_{\bm{\hat{n}}} \otimes \bm{\hat{n}} \right] : \bm{\tau} \cdot \left[ \bm{\hat{n}} \cdot \left( {\delta} \bm{u} \otimes {\bm{{\nabla}}}^{\alpha_1} \right) \right] \mathrm{d}\mathbb{A}
+ \oint_{\overline{\Gamma}} \left[ \left[ \bm{I}^{1-\alpha_1}_{\bm{\hat{m}}} \cdot ( \bm{I}^{1-\alpha_2}_{\bm{\hat{n}}} \cdot \bm{\tau} ) \right] \right] \cdot \delta \bm{u} \mathrm{d}\mathrm{l}
\end{split}
\end{equation}
The detailed derivation of the above governing equations is provided in the SI. In Eq.~(\ref{eq: PE_variation}), the tensor $\bm{R}$ is the projector onto the surface $\partial\Omega$, $\bm{\hat{m}}$ is the co-normal vector at the edges and [[$\cdot$]] operator denotes difference of the argument across both sides of the edge $\overline{\Gamma}$. For smooth edges (for example, a cube with filleted edges), the line integral vanishes analogous to classical formulations \cite{yurkov2011elastic}. $\bm{R}$ and $\bm{\hat{m}}$ are given as:
\begin{subequations}
    \label{eq: R_m}
    \begin{equation}
\bm{R}=\bm{1}-\bm{\hat{n}}\otimes\bm{\hat{n}}
\end{equation}
\begin{equation}
\bm{\hat{m}}=\bm{\hat{s}}\wedge\bm{\hat{n}}
\end{equation}
\end{subequations}
where $\bm{\hat{s}}$ is a unit vector tangent to the edge $\overline{\Gamma}$ and $\wedge$ denotes the exterior product. The operator $\bm{I}^{1-\alpha_m}_{\bm{\hat{n}}}(\cdot)$ is defined as:
\begin{equation}
\label{eq: integral_normal_operator}
    \bm{I}^{1-\alpha_m}_{\bm{\hat{n}}}(\cdot) = n_x {I}^{1-\alpha_m}_{x}(\cdot) \hat{x} + n_y {I}^{1-\alpha_m}_{y}(\cdot) \hat{y} + n_z {I}^{1-\alpha_m}_{z}(\cdot) \hat{z}
\end{equation}
such that $\bm{\hat{n}}=n_x \hat{x} + n_y \hat {y} + n_z \hat{z}$. The same definition directly extends to the operator $\bm{I}^{1-\alpha_m}_{\bm{\hat{m}}}(\cdot)$ that appears in Eq.~(\ref{eq: PE_variation}). Further, ${I}^{1-\alpha_m}_{x_j}(\cdot)$ is a Riesz integral operator defined in the following manner: 
\begin{equation}
\label{eq: reisz integral_def}
    I^{1-\alpha_m}_{x_j} \chi =\frac{1}{2}\Gamma(2-\alpha_m) \left[ l_{B_{j}}^{\alpha_m-1} \left({}_{{x_j}-l_{B_{j}}}I_{{x_j}}^{1-\alpha_m} \chi \right) - l_{A_{j}}^{\alpha_m-1} \left({}_{{x_j}}I_{{x_j}+l_{A_{j}}}^{1-\alpha_m} \chi \right) \right]
\end{equation}
where $x_1=x$, $x_2=y$ and $x_3=z$. ${}_{x_j-l_{B_{j}}}I_{{x_j}}^{1-\alpha_m}\chi$ and ${}_{{x_j}}I_{{x_j}+l_{A_{j}}}^{1-\alpha_m}\chi$ are the left and right Riesz integrals (in the $x_j$ direction) to the order $\alpha_m$ of an arbitrary function $\chi$. Further, the gradient operator denoted by $\bm{\widetilde{\nabla}}^{\alpha_m}(\cdot)$ is a Riesz Riemann-Liouville gradient (analogous to the RC gradient $\bm{{\nabla}}^{\alpha_m}(\cdot)$ in Eq.~(\ref{eq: RC_grad_operator})) containing Riesz Riemann-Liouville derivatives instead of RC derivatives. More specifically,
\begin{equation}
    \label{eq: RRL_grad_operator}
    \bm{\widetilde{\nabla}}^{\alpha_m} (\cdot) = \mathfrak{D}_{x}^{\alpha_m} (\cdot) \hat{x} + \mathfrak{D}_{y}^{\alpha_m} (\cdot) \hat{y} + \mathfrak{D}_{z}^{\alpha_m} (\cdot) \hat{z}
\end{equation}
where $\mathfrak{D}_{x_j}^{\alpha_m}(\cdot)$ is the Riesz Riemann-Liouville derivative of order $\alpha_m$ which is defined as:
\begin{equation}
    \label{eq: r_rl_frac_der_def}
    \mathfrak{D}^{\alpha_m}_{x_j} \chi = \frac{1}{2}\Gamma(2-\alpha_m) \left[ l_{B_{j}}^{\alpha_m-1} \left({}^{~~~~~RL}_{{x_j} - l_{B_{j}}} D^{\alpha_m}_{{x_j}} \chi \right) - l_{A_{j}}^{\alpha_m-1} \left( {}^{RL}_{x_j} D^{\alpha_m}_{{x_j} + l_{A_{j}}} \chi \right)\right]
\end{equation}
where ${}^{~~~~~RL}_{{x_j} - l_{B_{j}}} D^{\alpha_m}_{{x_j}} \chi$ and ${}^{RL}_{x_j} D^{\alpha_m}_{{x_j} + l_{A_{j}}} \chi$ are the left- and right-handed Riemann Liouville derivatives of $\chi$ to the order $\alpha_m$, in the $x_j$ direction. Note that the Riesz fractional derivative $\mathfrak{D}^{\alpha_m}_{x_j}(\cdot)$ and the Riesz fractional integral $I^{1-\alpha_m}_{x_j}(\cdot)$ are defined over the interval $({x_j}-l_{B_{j}},{x_j}+l_{A_{j}})$ unlike the RC fractional derivative $D^{\alpha_m}_{x_j}(\cdot)$ which is defined over the interval $({x_j}-l_{A_{j}},{x_j}+l_{B_{j}})$. This change in the terminals of the interval of the Riesz Riemann-Liouville integral and derivative follows from the variational simplifications (see SI). 

The first variation of the external work done follows directly from Eq.~(\ref{eq: 3D_ExtW}) as:
\begin{equation}
\label{eq: work_variation}
\delta \mathcal{V} = \int_{\Omega} ({\overline{\bm{b}}} \cdot \delta \bm{u} ) \mathrm{d} \mathbb{V} +\int_{\partial \Omega} ({\overline{\bm{t}}} \cdot \delta \bm{u} +\bm{\overline{q}} \cdot \bm{\hat{n}} \cdot ( \delta \bm{u} \otimes \bm{\nabla}^{\alpha_1})) \mathrm{d} \mathbb{A} + \oint_{\Gamma} (\bm{\overline{r}} \cdot \delta \bm{u}) \mathrm{d}\mathrm{l}
\end{equation}
Further, the first variation of the kinetic energy is obtained as:
\begin{equation}
    \label{eq: KE_variation}
    \delta T = - \int_{\Omega} \rho \bm{\ddot{u}} \cdot \delta \bm{u}\mathrm{d}\mathbb{V} 
\end{equation}
Now by using the extended Hamilton's principle in Eq.~(\ref{eq: extended_hamiltons_principle}) and applying the fundamental theorem of variational calculus, the elastodynamic governing equations for the 3D nonlocal continuum are obtained as:
\begin{equation}
\label{eq: GDE}
\bm{\widetilde{\nabla}}^{\alpha_1} \cdot \big( \bm{\sigma} - {\bm{\widetilde{\nabla}}}^{\alpha_2} \cdot \bm{\tau} \big) + \bm{\overline{b}} = \rho \bm{\ddot{u}} ~~~~ \forall ~ \bm{x} \in \Omega
\end{equation}
The associated boundary conditions are obtained as:
\begin{subequations}
\label{eq: BCs}
\begin{equation}
\bm{I}^{1-\alpha_1}_{\bm{\hat{n}}} \cdot \big( \bm{\sigma} - {\bm{\widetilde{\nabla}}}^{\alpha_2} \cdot \bm{\tau} \big) - \left[ \bm{R} : \left( \bm{I}^{1-\alpha_2}_{\bm{\hat{n}}} \cdot \bm{\tau} \right) \otimes {\bm{\widetilde{\nabla}}}^{\alpha_1} \right] = \bm{\overline{t}}  ~~~~ \text{or} ~~~~ {\bm{u}} = \bm{\overline{u}} ~~~~ \forall ~ \bm{x} \in \partial\Omega
\end{equation}
\begin{equation}
\left[ \bm{I}^{1-\alpha_2}_{\bm{\hat{n}}} \otimes \bm{n} \right] : \bm{\tau} = \overline{\bm{q}} ~~~~ \text{or} ~~~~ \bm{\hat{n}} \cdot ({\delta} \bm{u} \otimes \bm{\nabla}^{\alpha_1}) = \bm{\hat{n}} \cdot ({\delta} \bm{\overline{u}} \otimes \bm{\nabla}^{\alpha_1}) ~~~~ \forall ~ \bm{x} \in \partial\Omega
\end{equation}
\begin{equation}
\left[\left[ \bm{I}^{1-\alpha_1}_{\bm{\hat{m}}} \cdot ( \bm{I}^{1-\alpha_2}_{\bm{\hat{n}}} \cdot \bm{\tau} ) \right]\right]  = \bm{\overline{r}} ~~~~ \text{or} ~~~~ {\bm{u}} = \bm{\overline{u}} ~~~~ \forall ~ \bm{x} \in \overline{\Gamma}
\end{equation}
\end{subequations}
Note that the natural boundary conditions are nonlocal in nature. This is similar to what is seen in classical integral approaches \cite{eringen1972linear,reddy2014non}. The nonlocal nature follows from the nonlocal definition of the constitutive relations given in Eq.~(\ref{eq: constt}). It follows that the surface tractions depend on the response of a range of particles, hence leading to nonlocal boundary conditions. The partial horizon at the point $\bm{X}_3$ in Fig.~(\ref{fig: FCM})) serves as an example to illustrate the nonlocal nature of the boundary condition. We anticipate that the nonlocal nature of the natural boundary conditions does not concern us immediately as we will solve the above system of equations using a finite element (FE) technique. Recall that natural boundary conditions are implicitly satisfied when obtaining the solutions using FE techniques and are accurate up to the order of the specific finite element. Additionally, the following initial conditions are required to obtain the transient response:
\begin{equation}
\label{eq: IC}
    \delta \bm{u} = 0 ~~~~ \text{and} ~~~~ \delta \bm{\dot{u}}= 0 ~~~~ \forall~ \bm{x} \in \Omega~\text{at}~t=0
\end{equation}
Given the complex nature of the fractional-order governing equations and the associated boundary conditions, they do not generally admit closed-form analytical solutions. Consequently, numerical methods become indispensable to simulate the above governing equations. This issue is typical also of classical strain-gradient or integral nonlocal approaches, which typically are solved via numerical techniques \cite{askes2011gradient}.

In the following, we will use the fractional-order continuum formulation developed above to analyze both the static and the free vibration response of slender nonlocal structures, including a Timoshenko beam and a Mindlin plate. We will demonstrate that the fractional-order continuum model is able to capture both stiffening and softening effects depending on the values of the parameters involved in the fractional formulation. Numerical solutions will be obtained by using an adapted version of the fractional-order FEM (f-FEM) developed in \cite{patnaik2019FEM,patnaik2020plates} for fractional-order nonlocal BVPs. Note that the f-FEM is obtained by discretization of the Hamiltonian of the system using an isoparametric formulation. Hence, we only provide the weak form of the governing equations for the Timoshenko beam and the Mindlin plate.  The strong form of the governing equations for both the beam and the plate can be easily obtained following the detailed derivation of the 3D governing equations outlined in the SI.

%%%%%%%%%%%%%%%%%%%%%%%%%%%%%%%%%%%%%%%%%%%%%%%%%%%%%%%%%%%%%%%%%%%%%%%%%%%%%%%%%%%%%%%%%%%%%%%%%%%%%%
\section{Application to Timoshenko beams}
\label{sec: Timoshenko_beams}
\begin{figure}[h!]
    \centering
    \includegraphics[width=0.5\textwidth]{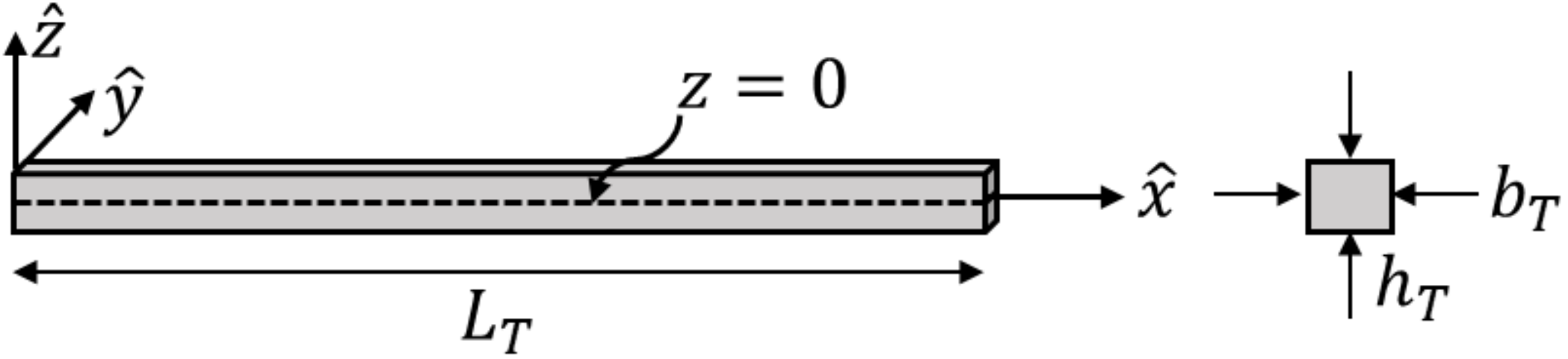}
    \caption{Schematic of the beam illustrating the different geometric parameters.}
    \label{fig: beam}
\end{figure}

We start analyzing the fractional-order continuum model by considering its application to a Timoshenko beam. A schematic of the undeformed beam along with the chosen Cartesian reference frame is illustrated in Fig.~(\ref{fig: beam}). The top surface of the plate is identified as $z=h_T/2$, while the bottom surface is identified as $z=-h_T/2$. The width of the beam is denoted as $b_T$. The domain corresponding to the mid-plane of the beam (i.e., $z=0$) is denoted as $\Omega_T$, such that $\Omega_T=[0,L_T]$ where $L_T$ is the length of the beam. The domain of the plate is identified by the tensor product $\Omega_T \otimes [-b_T/2,b_T/2] \otimes [-h_T/2,h_T/2]$. The subscript $T$ indicates that all the above dimensions correspond to the Timoshenko beam.

For the Timoshenko beam, analogous to the classical case, the axial and transverse components of the displacement field denoted by $u(x,y,z,t)$ and $w(x,y,z,t)$ at any spatial location $\bm{x}(x,y,z)$ are related to the mid-line displacements of the beam in the following manner:
\begin{subequations}
\label{eq: Timoshenko_kinematics}
    \begin{equation}
        u(x,y,z,t)=u_0(x,t)-z\theta_0(x,t)
    \end{equation}
    \begin{equation}
        w(x,y,z,t)=w_0(x,t)
    \end{equation}
\end{subequations}
where $u_0$ and $w_0$ are the mid-plane axial and transverse displacements of the beam, and $\theta_0$ is the rotation of the transverse normal of the beam about the $\hat{y}$ axis. In the following, for a compact notation, the functional dependence of the displacement fields on the spatial and the temporal variables will be implied unless explicitly expressed to be constant. Based on the above displacement fields, the non-zero strain components in the Timoshenko beam are evaluated using Eq.~(\ref{eq: infinitesimal_fractional_strain}) as:
\begin{subequations}
\label{eq: Timoshenko_strains}
\begin{equation}
    \epsilon_{xx}=D^{\alpha_1}_x u_0-zD^{\alpha_1}_x \theta_0
\end{equation}
\begin{equation}
    \epsilon_{xz}=\frac{1}{2}\left[D^{\alpha_1}_x w_0-\theta_0\right]
\end{equation}
\end{subequations}
The strain-gradients developed in the beam are obtained using Eq.~(\ref{eq: frac_straingrad}) as:
\begin{subequations}
\label{eq: Timoshenko_strain_gradients}
\begin{equation}
    \eta_{xxr}=D^{\alpha_2}_r\left[D^{\alpha_1}_x u_0-zD^{\alpha_1}_x \theta_0\right]
\end{equation}
\begin{equation}
    \eta_{xzr}=D^{\alpha_2}_r\left[\frac{1}{2}\left[D^{\alpha_1}_x w_0-\theta_0\right]\right]
\end{equation}
\end{subequations}
where $r\in\{x,y,z\}$. Specializing the above expressions, the following strain-gradient components are obtained exactly:
\begin{subequations}
\label{eq: Timoshenko_exact_strain_grads}
\begin{equation}
    \eta_{xxz} = D^{\alpha_2}_z [D^{\alpha_1}_x u_0] - D^{\alpha_2}_z[zD^{\alpha_1}_x \theta_0] = -D^{\alpha_1}_x \theta_0
\end{equation}
\begin{equation}
    \eta_{xzx} = \frac{1}{2} \left[ D^{\alpha_2}_x [D^{\alpha_1}_x w_0] - D^{\alpha_2}_x \theta_0 \right] = -\frac{1}{2} D^{\alpha_2}_x \theta_0
\end{equation}
\begin{equation}
    \eta_{xzz} = \frac{1}{2} \left[ D^{\alpha_2}_z [D^{\alpha_1}_x w_0] - D^{\alpha_2}_z \theta_0 \right] = 0
\end{equation}
\end{subequations}

In the above simplification we have used that $D^{\alpha_2}_z[z]=1$, which is exact and follows immediately from the definition of the RC derivative defined in Eq.~(\ref{eq: RC_definition}). Further, assuming small displacement gradients $(O(\varepsilon))$, the strain-gradient $\eta_{xxx}$ is $O(\varepsilon^2)$ while the strain-gradients in Eq.~(\ref{eq: Timoshenko_exact_strain_grads}) are either $O(\varepsilon)$ or exactly zero. Hence it appears that, for the normal strain $\epsilon_{xx}$, the transverse strain-gradient $\eta_{xxz}$ is significant when compared to the axial gradient $\eta_{xxx}$. Conversely, for the shear strain $\epsilon_{xz}$, the axial gradient $\eta_{xzx}$ is significant while the transverse gradient $\eta_{xzz}$ is exactly zero. Thus, when obtaining the response of the beam via the weak form, the contribution of the strain-gradient $\eta_{xxx}$ can be ignored when compared to the contribution of the non-zero strain-gradients in Eq.~(\ref{eq: Timoshenko_exact_strain_grads}). We have further justified this approximation in detail in the SI.

The first variations of the nonlocal potential energy, the work done by externally applied forces, and the kinetic energy corresponding to the Timoshenko beam assumptions are obtained as:
\begin{subequations}
\label{eq: Timo_virtual_quantities}
\begin{equation}
\label{eq: Timo_virtual_strain_energy}
     \delta \mathcal{U} = \int_0^{L_T} \Big[ N_{xx}\delta D^{\alpha_1}_x u_0 + M_{xx}\delta D^{\alpha_1}_x \theta_0 +  Q_{xz}\delta(D^{\alpha_1}_x w_0-\theta_0) + \overline{N}_{xxz} \delta D^{\alpha_1}_x \theta_0 + \overline{N}_{xzx} \delta D^{\alpha_2}_x \theta_0 \Big] \mathrm{d} x
\end{equation}
\begin{equation}
\label{eq: Timo_virtual_work}
    \delta \mathcal{V} = \int_0^{L_T} \int_{-\frac{b_T}{2}}^{\frac{b_T}{2}} \int_{-\frac{h_T}{2}}^{\frac{h_T}{2}}  \left[ F_x\delta u_0 + F_z \delta w_0 + M_{\theta_0} \delta \theta_0 \right]\mathrm{d}z \mathrm{d} y \mathrm{d} x
\end{equation}
\begin{equation}
\label{eq: Timo_virtual_kinetic_energy}
    \delta T = \int_0^{L_T} \int_{-\frac{b_T}{2}}^{\frac{b_T}{2}} \int_{-\frac{h_T}{2}}^{\frac{h_T}{2}} \rho \left[ \big(\dot{u}_0 - z\dot{\theta_0} \big)\big(\delta\dot{u}_0 - z\delta\dot{\theta_0} \big) +\dot{w}_0 \delta \dot{w}_0 \right] \mathrm{d}z \mathrm{d} y \mathrm{d} x
\end{equation}
\end{subequations}
$\{F_x,F_z\}$ are the external loads applied in the axial ($\hat{x}$) and transverse ($\hat{z}$) directions, respectively, and $M_{\theta_0}$ is the external moment applied about the $\hat{y}$ axis.
The axial stress resultant $\{ N_{xx}\}$, the shear resultant $\{Q_{xz}\}$, the moment resultant $\{ M_{xx}\}$, and the higher-order stress resultants $\{ \overline{N}_{xxz}, \overline{N}_{xzx}\}$ in Eq.~(\ref{eq: Timo_virtual_strain_energy}) are given as:
\begin{equation}
\label{eq: Timo_stress_resultants}
    \{ N_{xx}, Q_{xz}, M_{xx}, \overline{N}_{xxz}, \overline{N}_{xzx}\} = \int_{-\frac{b_T}{2}}^{\frac{b_T}{2}} \int_{-\frac{h_T}{2}}^{\frac{h_T}{2}} \{ \sigma_{xx}, K_s \sigma_{xz}, -z\sigma_{xx}, -\tau_{xxz}, -K_s \tau_{xzx}\} \mathrm{d}z \mathrm{d} y
\end{equation}
where $K_s$ is the shear correction factor.

In the following, we briefly discuss the f-FEM method used to numerically simulate the fractional-order system. The details of the f-FEM are extensive and will not be reported here, but the interested reader can refer to \cite{patnaik2019FEM,patnaik2020plates}. The f-FEM for the Timoshenko beam is formulated by obtaining a discretized form of the first variation of the Lagrangian of the beam. For this purpose, the beam domain $\Omega_T=[0,L]$ is uniformly discretized into disjoint three-noded line elements and the different fractional derivatives that appear in Eq.~(\ref{eq: Timo_virtual_strain_energy}) are expressed as:
\begin{subequations}
\label{eq: FE_frac_derivatives_final}
\begin{equation}
    D^{\alpha_1}_{x}\left[u_0(x)\right] = [\tilde{B}^{\alpha_1}_{u_0,x}(x)]\{U\}
\end{equation}
\begin{equation}
    D^{\alpha_m}_{x}\left[\theta_0(x)\right] = [\tilde{B}^{\alpha_m}_{\theta_0,x}(x)] \{U\}
\end{equation}
\begin{equation}
   D^{\alpha_1}_{x}\left[w_0(x)\right] - \theta_x = \left[ [\tilde{B}^{\alpha_1}_{w_0,x}(x)] - [{\mathbb{L}}^{(\theta_0)}(x)] \right] \{U\}
\end{equation}
\end{subequations}
where $\{U\}$ denotes the global degrees of freedom vector and $[{\mathbb{L}}^{(\theta_0)}(x)]$ is obtained by assembling the element interpolation vectors for $\theta_0$. The matrices $[\tilde{B}^{\alpha_m}_{\square,x}(x)]$ contain the fractional-order derivatives of the shape functions used to interpolate the nodal displacement degrees of freedom of the Timoshenko beam. A brief discussion on the details of these matrices is provided in SI. By using the above expressions for the FE approximation of the different fractional-order derivatives, the first variation of the potential energy $\delta\mathcal{U}$ given in Eq.~(\ref{eq: Timo_virtual_strain_energy}) is obtained as:
\begin{equation}
\label{eq: FE_expression_virtual_strain}
    \delta \mathcal{U} = \delta \{U\}^T \left[ \int_{\Omega_T} [\tilde{B}_T(x)]^T [S_T] [\tilde{B}_T(x)] \mathrm{d}\Omega_T \right] \{U\} = \delta \{U\}^T [K_T] \{U\} 
\end{equation}
where $[S_T]$ is the constitutive matrix of the beam and $[\tilde{B}_T(x)]$ is given as:
\begin{equation}
\label{eq: FE_B_matrices_collection}
\begin{split}
        [\tilde{B}_T(x)] = \bigg[ \underbrace{[\tilde{B}^{\alpha_1}_{u_0,x}(x)]^T, \left[ [\tilde{B}^{\alpha_1}_{w_0,x}(x)] - [{\mathbb{L}}^{(\theta_0)}(x)] \right]^T}_{\text{Contributions from the nonlocal strains}}, \underbrace{\vphantom{[\tilde{B}^{\alpha_1}_{u_0,x}(x)]^T, \left[ [\tilde{B}^{\alpha_1}_{w_0,x}(x)] - [{\mathbb{L}}^{(\theta_0)}(x)] \right]^T}[\tilde{B}^{\alpha_1}_{\theta_0,x}(x)]^T, [\tilde{B}^{\alpha_1}_{\theta_0,x}(x)]^T, [\tilde{B}^{\alpha_2}_{\theta_0,x}(x)]^T}_{\text{Contributions from the nonlocal strain-gradients}} \bigg]^T 
\end{split}
\end{equation}
The algebraic equations for the f-FEM are given as:
\begin{equation}
    \label{eq: FE_algebraic_equations}
    [M_T]\{\ddot{U}\} + [K_T]\{U\} = \{F_T\}
\end{equation}
where the stiffness matrix $[K_T]$ is indicated in Eq.~(\ref{eq: FE_expression_virtual_strain}). The expressions for the force vector $\{F_T\}$ and the mass matrix $[M_T]$ follow directly from classical Timoshenko beam formulations and are provided in SI. The solution of the algebraic Eq.~\eqref{eq: FE_algebraic_equations} gives the nodal displacement variables which can then be used along with the kinematic relations in Eq.~(\ref{eq: Timoshenko_kinematics}) to determine the displacement field at any point in the beam. Note that the f-FEM also involves the numerical evaluation of the mass matrix, the stiffness matrix, and the force vector. The procedure to numerically evaluate the mass matrix and the force vector follows directly from classical FE formulations. The stiffness matrix of the fractional-order nonlocal system requires the evaluation of the different nonlocal matrices given in Eq.~(\ref{eq: FE_B_matrices_collection}). Further, the attenuation function in the fractional-order model involves an end-point singularity due to the nature of the kernel \cite{kilbas2006theory}. The fractional-order nonlocal interactions as well as the end-point singularity are addressed in detail in \cite{patnaik2019FEM,patnaik2020plates}. We emphasize that the numerical integration procedure presented in \cite{patnaik2019FEM,patnaik2020plates} directly extends to the evaluation of the stiffness matrix of the FE governing equations derived in this study.

%%%%%%%%%%%%%%%%%%%%%%%%%%%%%%%%%%%%%%%%%%%%%%%%%%%%%%%%%%%%%%%%%%%%%%%%%%%%%%%%%%%%%%%%%%%%%%%%%%%%%%
\subsection{Static response}
\label{ssec: Timoshenko_static}

In this section, we analyse the static response of the Timoshenko beam which was obtained by solving the static part of the fractional-order FE algebraic equations in Eq.~(\ref{eq: FE_algebraic_equations}). In the following study, the dimensions of the beam were fixed to be $L_T=1$m, $b_T=0.1$m and $h_T=0.05$m $(=L/20)$. The simplified constitutive relations proposed in \cite{lazar2006dislocations} were used in this study:
\begin{equation}
    \label{eq: simplfied_CE}
    \mathcal{U} = \frac{1}{2} C_{ijkl} \epsilon_{ij} \epsilon_{kl} + \frac{1}{2} l^2_* C_{ijmn} \eta_{mnk} \eta_{ijk}
\end{equation}
The material was assumed to be isotropic with an elastic modulus $E=30$GPa, Poisson's ratio $\nu=0.3$ and density $\rho=2700$kg/m${}^3$. Further, we have assumed a symmetric and isotropic horizon of nonlocality for points sufficiently inside the domain of the beam, that is $l_{A_x} = l_{B_x} = l_f$. For points located close to a boundary, the length scales are truncated as shown in Fig.~(\ref{fig: FCM}). Using the above material properties, we analyzed the effect of the following fractional model parameters: nonlocal strain order ($\alpha_1$), strain-gradient order ($\alpha_2$), nonlocal horizon length ($l_f$) and microstructure length ($l_*$), on the static response of the Timoshenko beam. We merely note that the Young’s modulus $E$ and the Poisson’s ratio $\nu$ chosen above correspond to a general class of soft metals (e.g. lead). Given the linearity of the problem and the fact that results will be presented in a normalized form, the choice of specific elastic constants is quite immaterial for the interpretation of the results.

We analyzed the static response of the beam subject to a uniformly distributed transverse load (UDTL) of magnitude $F_z=10^{7}\text{N}/\text{m}$ for two different kinds of boundary conditions: 1) clamped-clamped (CC), and 2) simply supported at both ends (SS). For each boundary condition, we obtained the response of the beam for the following different cases:
\begin{itemize}
    \item Case 1: the fractional-orders $\alpha_1$ and $\alpha_2$ were varied within the range [0.7,1] for fixed values of the nonlocal horizon length $l_f=0.5\text{m}~(=L/2)$. For this case, the microstructural length was chosen as $l_*=0.005\text{m}~(L/200)$ for the CC beam and $l_*=0.002\text{m}~(L/500)$ for the SS beam.
    
    \item Case 2: the horizon length $l_f$ was varied within the range $[0.1,0.5]\text{m}~ (=[L/10,L/2])$ and the microstructure length $l_*$ was varied in $[0.002,0.01]\text{m}~(=[L/500,L/100])$, for fixed values of the fractional-orders $\alpha_1=\alpha_2=0.8$. Both orders were chosen in the fractional range so to obtain more general conditions (see Fig.~(\ref{fig: parameter_map})).
\end{itemize}
We emphasize that, while the choice of the different fractional-model parameters were somewhat arbitrary, their specific value does not affect the generality of the results. The range of the fractional-orders $\alpha_1$ and $\alpha_2$ was selected following the restriction in Eq.~(\ref{eq: restriction_on_order}). The specific ranges for $l_f$ and $l_*$ were chosen in order to demonstrate the ability of the fractional-order framework in capturing both stiffening and softening effects.

\begin{figure}[h!]
    \centering
    \includegraphics[width=\textwidth]{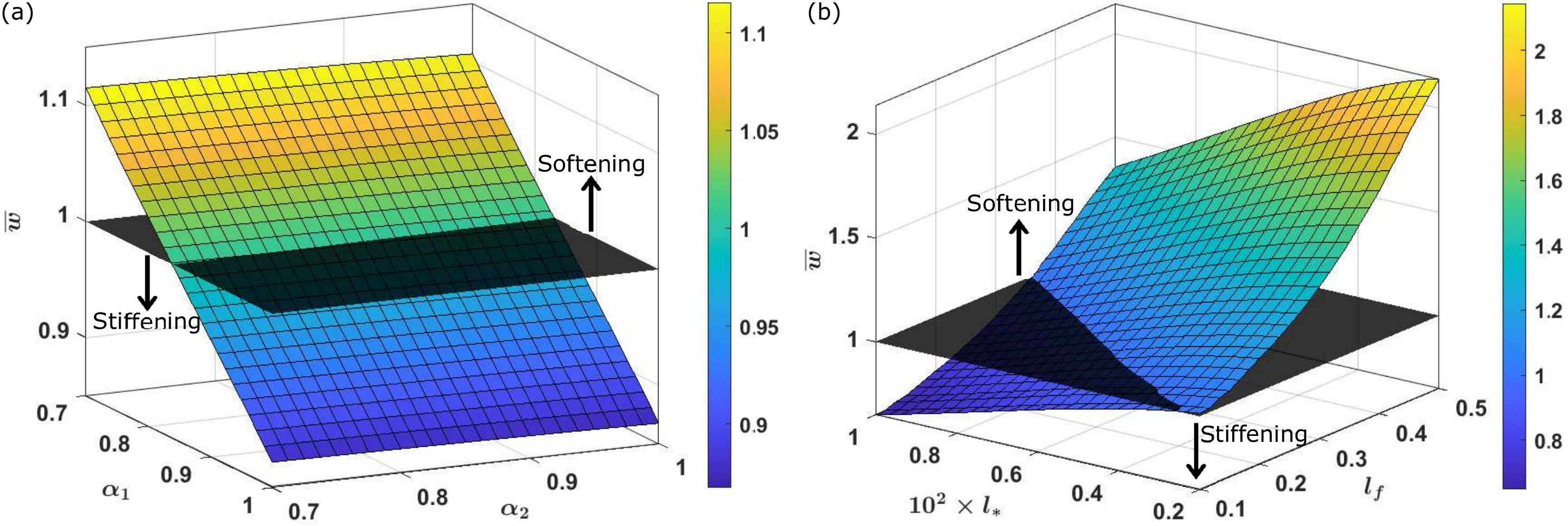}
    \caption{Non-dimensional transverse displacement at the center point of the Timoshenko beam subject to clamped-clamped boundary conditions. Results are obtained via the fractional-order formulation. The response is parameterized for different values of (a) the fractional orders and (b) the length scales.}
    \label{fig: Timoshenko_CC}
\end{figure}

\begin{figure}[h!]
    \centering
    \includegraphics[width=\textwidth]{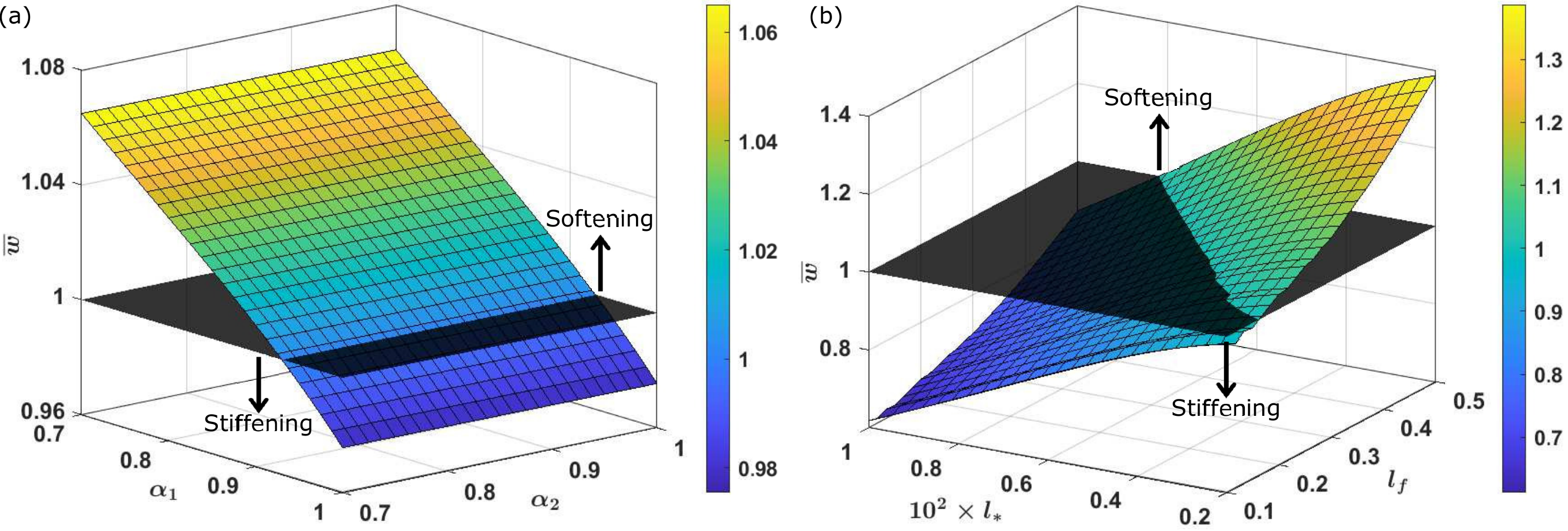}
    \caption{Non-dimensional transverse displacement at the center point of the Timoshenko beam subject to simply supported boundary conditions. Results are obtained via the fractional-order formulation and parameterized for different values of (a) the fractional orders and (b) the length scales.}
    \label{fig: Timoshenko_SS}
\end{figure}

The numerical results, expressed in terms of the static transverse displacement and corresponding to Case 1 for the CC beam and the SS beam, are presented in Fig.~(\ref{fig: Timoshenko_CC}a) and Fig.~(\ref{fig: Timoshenko_SS}a), respectively. Similarly, the results for Case 2 subject to either CC or SS boundary conditions are provided in Fig.~(\ref{fig: Timoshenko_CC}b) and Fig.~(\ref{fig: Timoshenko_SS}b), respectively. The results presented for each case correspond to the maximum transverse displacement observed in the beam at the mid point ($w_0(L_T/2)$). To clearly visualize the extent of softening and stiffening occurring in the beam, the maximum transverse displacement was non-dimensionalized against the maximum transverse displacement obtained for a classical Timoshenko beam in the absence of both nonlocal and strain-gradient effects. More specifically, the non-dimensional transverse displacement $(\overline{w})$ for each specific boundary configuration, was obtained by dividing the maximum transverse displacement of the fractional-order beam by the maximum transverse displacement of the classical beam for the same boundary condition. The maximum transverse displacement obtained for the classical CC beam was ${w}_0=8.95\times10^{-2}\text{m}$ and for the classical SS beam was ${w}_0=41.93\times10^{-2}\text{m}$. Note that a higher value of the static displacement with respect to the classical solution indicates softening of the structure, while a lower value of the transverse displacement indicates a stiffening of the structure.

As evident from Fig.~(\ref{fig: Timoshenko_CC}) and  Fig.~(\ref{fig: Timoshenko_SS}), the fractional-order continuum formulation is able to capture both stiffening and softening response of the Timoshenko beam depending on the choice of the nonlocal parameters. Note that the horizontal reference plane in black color denotes the non-dimensional classical solution ($\overline{w}=1$). When the transverse displacement is above this plane ($\overline{w}>1$) it indicates a softened response while, values below the plane ($\overline{w}<1$) indicate a stiffened response. The results presented for the different cases lead to the following conclusions on the specific effects of the different fractional model parameters:
\begin{itemize}
    \item \textbf{Effect of $\alpha_1$:} As discussed in \cite{patnaik2019FEM}, a decrease in the value of $\alpha_1$ leads to an increase in the strength of the power-law kernel that captures nonlocal interactions across the horizon of nonlocality. Consequently, the resulting formulation exhibits a greater degree of softening with respect to the classical response. Recall that for $\alpha_1=1$ and $l_*=0$ (no microstructural effects), the classical local continuum formulation is recovered from the fractional-order formulation.
    \item \textbf{Effect of $l_f$:} recall that $l_f$ indicates the size of the nonlocal horizon, thus by increasing the value of $l_f$ the size of the horizon of nonlocality increases. It follows that a larger number of points within the solid is accounted contribute to the nonlocal interactions, thus the degree of nonlocality increases and so does the degree of softening of the structure.
    \item \textbf{Effect of $\alpha_2$:} Recall from \S\ref{ssec: 1D_model} that the strain-gradient order $\alpha_2$ captures the nonlocal effects of the strain-gradients. Thus, analogous to $\alpha_1$, a decrease in the value of $\alpha_2$ leads to an increase in the strength of the power-law kernel that captures nonlocal strain-gradient contributions across the horizon of nonlocality. Consequently, the resulting formulation would exhibit a softening with respect to the classical first-order strain gradient response. Note that for $\alpha_1=1$ and $\alpha_2=1$, the classical first-order strain-gradient theory is recovered from the fractional-order formulation.
    \item \textbf{Effect of $l_*$:} As evident from the discussion of the lattice structure in \S\ref{ssec: 1D_model}, the microstructural length parameter $l_*$ plays the same role as in classical strain-gradient formulations. Thus, an increase in the value of $l_*$ leads to a stiffer response of the structure. 
\end{itemize}
The effects discussed above are schematically summarized in Fig.~(\ref{fig: parameter_map}), which provides a visual representation of the resulting formulation as a function of the different parameters.

\begin{figure}[h!]
    \centering
    \includegraphics[width=\textwidth]{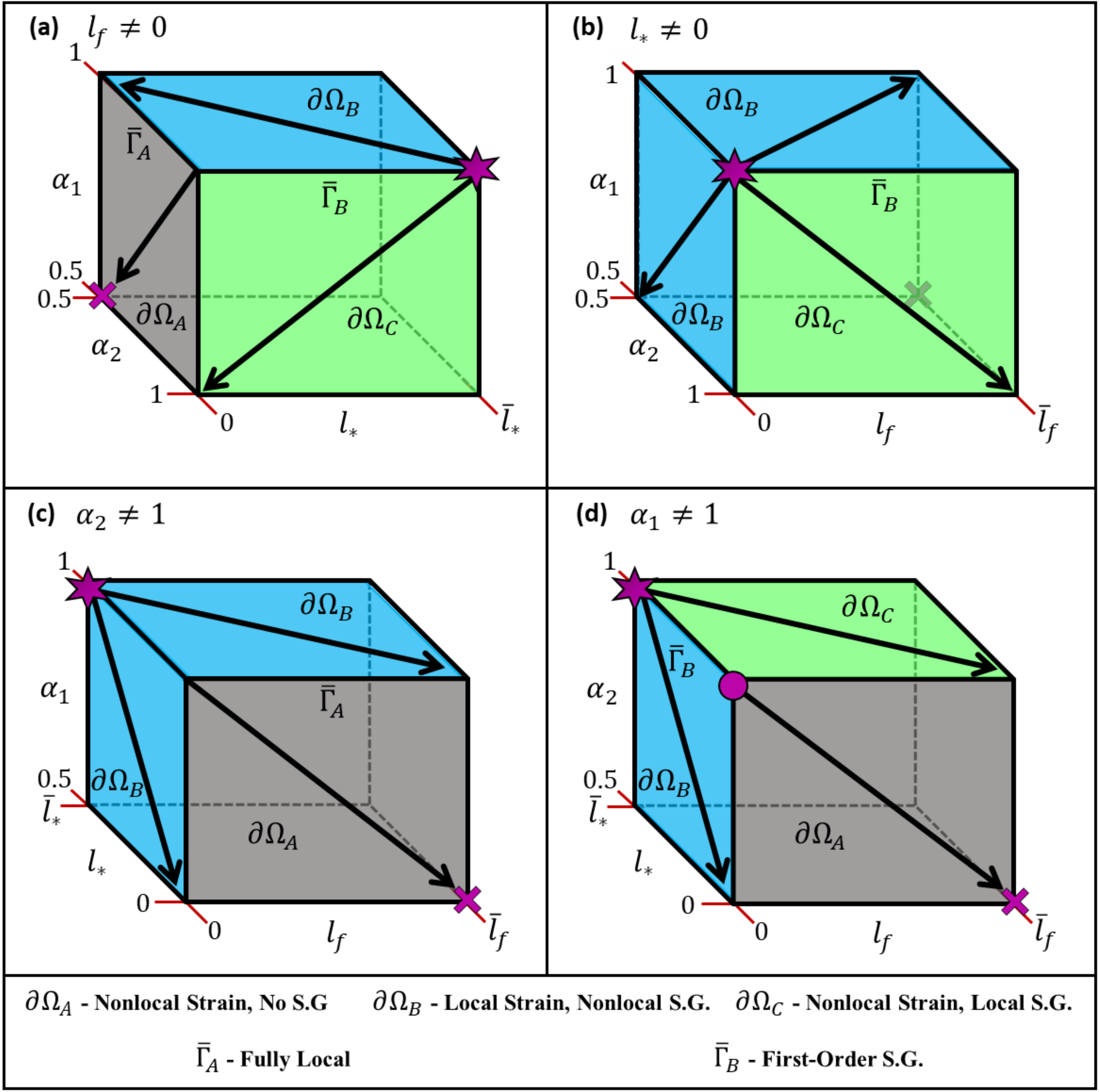}
    \caption{Schematic illustrating the effects of the different fractional-model parameters on the direction of softening or stiffening. In the above figure, S.G. denotes strain gradient, and $\overline{l}_f$ and $\overline{l}_*$ indicate the upper bound on the nonlocal horizon length and the microstructural length. The direction of the solid arrow lying on a particular plane, indicates the direction of softening. It is immediate that the opposite direction would lead to stiffening. In each sub-figure, the combination of parameters that would result in the stiffest and the softest solution is indicated by a six-edged star symbol ($\ast$) and a cross symbol ($\times$), respectively. In (d) the fully local solution is obtained at the corner indicated by a filled circular symbol.}
    \label{fig: parameter_map}
\end{figure}

%%%%%%%%%%%%%%%%%%%%%%%%%%%%%%%%%%%%%%%%%%%%%%%%%%%%%%%%%%%%%%%%%%%%%%%%%%%%%%%%%%%%%%%%%%%%%%%%%%%%%%
\subsection{Free vibration response}
\label{ssec: Timoshenko_free_vibration}
\begin{figure}[h!]
    \centering
    \includegraphics[width=\textwidth]{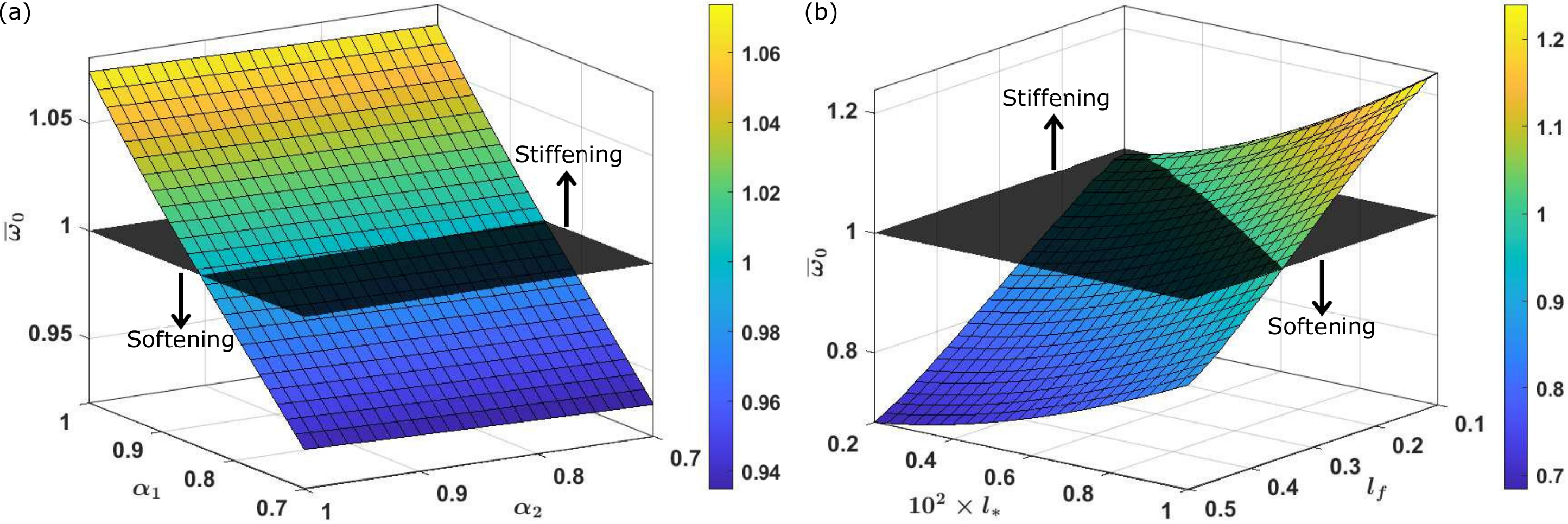}
    \caption{Non-dimensional natural frequency of the Timoshenko beam clamped at its boundaries and parameterized for different values of (a) the fractional orders and (b) the length scales.}
    \label{fig: Eigen_Timoshenko_CC}
\end{figure}

\begin{figure}[h!]
    \centering
    \includegraphics[width=\textwidth]{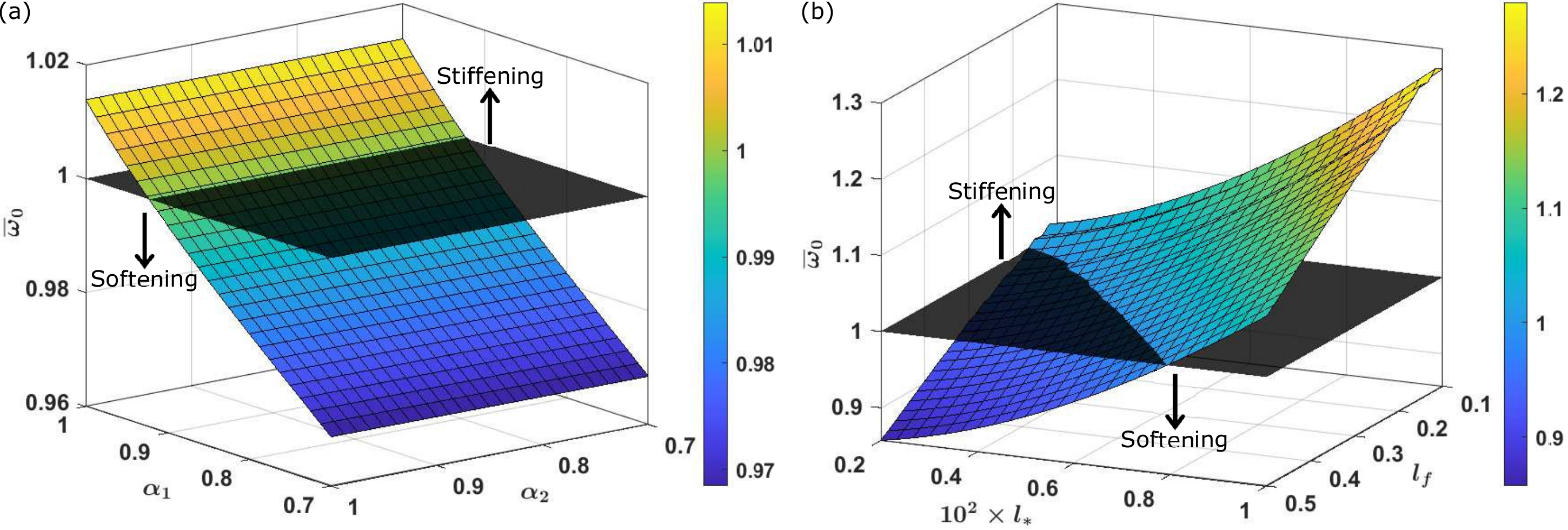}
    \caption{Non-dimensional natural frequency of the Timoshenko beam simply-supported at its boundaries and parameterized for different values of (a) the fractional orders and (b) the length scales.}
    \label{fig: Eigen_Timoshenko_SS}
\end{figure}

In the interest of a comprehensive analysis, we analyse the effect of the different fractional model parameters on the natural frequency of transverse vibration of the Timoshenko beam. The material properties chosen for this study are the same as those provided for the static study in \S\ref{ssec: Timoshenko_static}. The natural frequencies are obtained by solving the eigenvalue problem:
\begin{equation}
    \label{eq: eigen_system}
    [M_T]^{-1}[K_T]\{\mathbb{U}\}=\omega_0^2\{\mathbb{U}\}
\end{equation}
which is derived by assuming a periodic solution $\{U\}=\{\mathbb{U}\}e^{-i\omega_0t}$ to the homogeneous part of the algebraic FE Eq.~(\ref{eq: FE_algebraic_equations}). In the above assumed solution, $\omega_0$ denotes the natural frequency of vibration, and $\mathbb{U}$ is the amplitude of the harmonic oscillation. Similar to \S\ref{ssec: Timoshenko_static}, we obtained the natural frequencies of CC and SS beams for the two different cases: Case 1 and Case 2. The results are presented in Figs.~(\ref{fig: Eigen_Timoshenko_CC}, \ref{fig: Eigen_Timoshenko_SS}). Similar to the static analysis, the natural frequency obtained for each case ($\overline{\omega}_0$) was non-dimensionalized against the natural frequency of a classical local beam, which was found to be $53$Hz for the CC beam and $24$Hz for the SS beam. Note that a lower value of the natural frequency ($\overline{\omega}_0<1$) with respect to the classical solution indicates softening of the structure, while a higher value of the natural frequency ($\overline{\omega}_0<1$) indicates a stiffening of the structure. Clearly, the results presented in Figs.~(\ref{fig: Eigen_Timoshenko_CC}, \ref{fig: Eigen_Timoshenko_SS}) complement the discussion presented in \S\ref{ssec: Timoshenko_static}, on the effect of the different fractional-order parameters on the static response of the beam.

%%%%%%%%%%%%%%%%%%%%%%%%%%%%%%%%%%%%%%%%%%%%%%%%%%%%%%%%%%%%%%%%%%%%%%%%%%%%%%%%%%%%%%%%%%%%%%%%%%%%%%
\section{Application to Mindlin plates}
\label{sec: Minldin_plates}
\begin{figure}[h]
    \centering
    \includegraphics[width=0.5\textwidth]{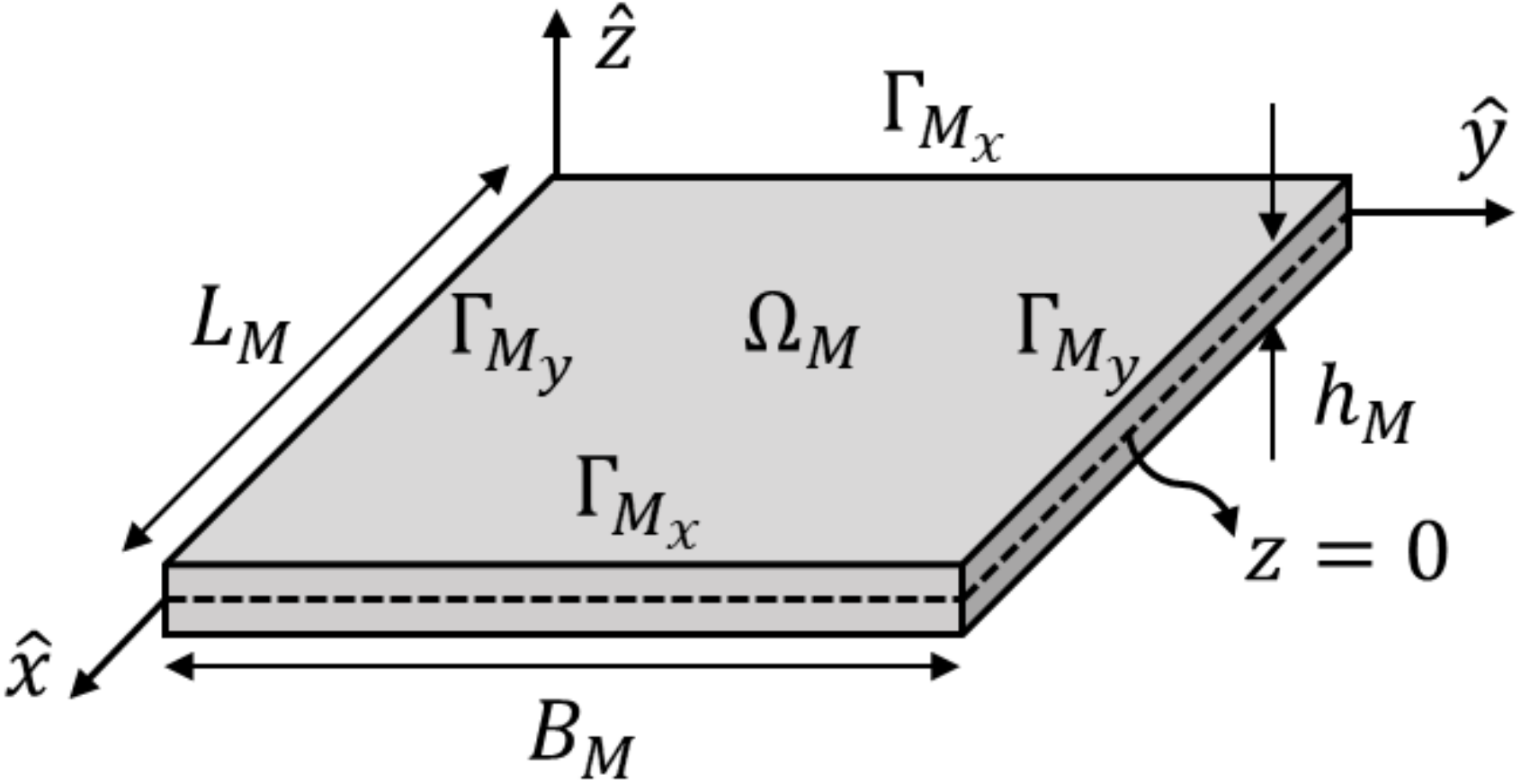}
    \caption{Schematic of the rectangular plate illustrating the different geometric parameters.}
    \label{fig: plate}
\end{figure}

We extend the studies carried out in \S\ref{sec: 3D continuum} and \S\ref{sec: Timoshenko_beams} to develop a fractional-order analogue of the classical Mindlin plate formulation that captures both stiffening and softening response. A schematic of the undeformed rectangular plate along with the chosen Cartesian reference frame is given in Fig.~(\ref{fig: plate}). The top surface of the plate is identified as $z=h_M/2$, while the bottom surface is identified as $z=-h_M/2$. The domain corresponding to the mid-plane of the plate (i.e., $z=0$) is denoted as $\Omega_M$, such that $\Omega_M=[0,L_M]\otimes[0,B_M]$ where $L_M$ and $B_M$ are the length and width of the plate, respectively. The domain of the plate is identified by the tensor product $\Omega_M \otimes [-h_M/2,h_M/2]$. The edges forming the boundary of the mid-plane of the plate are denoted as $\{\Gamma_{M_x},\Gamma_{M_y}\}$. The subscript $M$ indicates that all the above dimensions correspond to the Mindlin plate.

For the Mindlin plate, following the coordinate system illustrated in Fig.~(\ref{fig: plate}), the in-plane and transverse components of the displacement field, denoted by $u(x,y,z,t)$, $v(x,y,z,t)$ and $w(x,y,z,t)$ at any spatial location $\bm{x}(x,y,z)$, are related to the mid-plane displacements of the plate in the following manner:
\begin{subequations}
\label{eq: Mindlin_Kinematics}
    \begin{equation}
    u(x,y,z,t)=u_0(x,y,t) - z\theta_x(x,y,t)
    \end{equation}
    \begin{equation}
    v(x,y,z,t)=v_0(x,y,t) - z\theta_y(x,y,t)
    \end{equation}
    \begin{equation}
    w(x,y,z,t)=w_0(x,y,t)
    \end{equation}
\end{subequations}
where $u_0$, $v_0$, and $w_0$ are the mid-plane displacements of the plate along the $\hat{x}$, $\hat{y}$, and $\hat{z}$ directions. $\theta_x$ and $\theta_y$ are the rotations of the transverse normal about the $\hat{y}$ and $\hat{x}$ axes, respectively. In the interest of a more compact notation, the functional dependence of the displacement fields on the spatial and the temporal variables will be implied unless explicitly expressed to be constant. Based on the above displacement fields, the non-zero strain components in the fractional-order Mindlin plate are evaluated using Eq.~(\ref{eq: infinitesimal_fractional_strain}) as:
\begin{subequations}
\label{eq: Mindlin_strains}
\begin{equation}
    \epsilon_{xx}=D^{\alpha_1}_x u_0-zD^{\alpha_1}_x \theta_x
\end{equation}
\begin{equation}
    \epsilon_{yy}=D^{\alpha_1}_y v_0-zD^{\alpha_1}_y \theta_y
\end{equation}
\begin{equation}
    \epsilon_{xy}=\frac{1}{2}\left[D^{\alpha_1}_y u_0+D^{\alpha_1}_x v_0-z(D^{\alpha_1}_y \theta_x+D^{\alpha_1}_x \theta_y)\right]
\end{equation}
\begin{equation}
    \epsilon_{xz}=\frac{1}{2}\left[D^{\alpha_1}_x w_0-\theta_x\right]
\end{equation}
\begin{equation}
    \epsilon_{yz}=\frac{1}{2}\left[D^{\alpha_1}_y w_0-\theta_y\right]
\end{equation}
\end{subequations}
The strain-gradients developed in the plate are obtained using Eq.~(\ref{eq: frac_straingrad}) as:
\begin{subequations}
\label{eq: Mindlin_strain_gradients}
\begin{equation}
    \eta_{xxr}=D^{\alpha_2}_r\left[D^{\alpha_1}_x u_0-zD^{\alpha_1}_x \theta_x\right]
\end{equation}
\begin{equation}
    \eta_{yyr}=D^{\alpha_2}_r\left[D^{\alpha_1}_y v_0-zD^{\alpha_1}_y \theta_y\right]
\end{equation}
\begin{equation}
    \eta_{xyr}=D^{\alpha_2}_r\left[\frac{1}{2}\left[D^{\alpha_1}_y u_0+D^{\alpha_1}_x v_0-z(D^{\alpha_1}_y \theta_x+D^{\alpha_1}_x \theta_y)\right]\right]
\end{equation}
\begin{equation}
    \eta_{xzr}=D^{\alpha_2}_r\left[\frac{1}{2}\left[D^{\alpha_1}_x w_0-\theta_x\right]\right]
\end{equation}
\begin{equation}
    \eta_{yzr}=D^{\alpha_2}_r\left[\frac{1}{2}\left[D^{\alpha_1}_y w_0-\theta_y\right]\right]
\end{equation}
\end{subequations}
where $r\in\{x,y,z\}$. While simplifying the expressions in the above equation, the following strain-gradients are obtained exactly:
\begin{subequations}
\label{eq: mindlin_exact_strain_grads}
\begin{equation}
    \eta_{xxz} = D^{\alpha_2}_z [D^{\alpha_1}_x u_0] - D^{\alpha_2}_z[zD^{\alpha_1}_x \theta_x] = -D^{\alpha_1}_x \theta_x
\end{equation}
\begin{equation}
    \eta_{yyz} = D^{\alpha_2}_z [D^{\alpha_1}_y v_0] -D^{\alpha_2}_z[zD^{\alpha_1}_y \theta_y] = -D^{\alpha_1}_y \theta_y
\end{equation}
\begin{equation}
    \eta_{xyz} = \frac{1}{2}\left[ D^{\alpha_2}_z [D^{\alpha_1}_y u_0 + D^{\alpha_1}_x v_0] - D^{\alpha_2}_z [zD^{\alpha_1}_y \theta_x + zD^{\alpha_1}_x \theta_y]\right] = -\frac{1}{2}[D^{\alpha_1}_y \theta_x + D^{\alpha_1}_x \theta_y]
\end{equation}
\begin{equation}
    \eta_{xzx} = \frac{1}{2} \left[ D^{\alpha_2}_x [D^{\alpha_1}_x w_0] - D^{\alpha_2}_x \theta_x \right] = -\frac{1}{2} D^{\alpha_2}_x \theta_x
\end{equation}
\begin{equation}
    \eta_{xzy} = \frac{1}{2} \left[ D^{\alpha_2}_y [D^{\alpha_1}_x w_0] - D^{\alpha_2}_y \theta_x \right] = -\frac{1}{2} D^{\alpha_2}_y \theta_x
\end{equation}
\begin{equation}
    \eta_{xzz} = \frac{1}{2} \left[ D^{\alpha_2}_z [D^{\alpha_1}_x w_0] - D^{\alpha_2}_z \theta_x \right] = 0
\end{equation}
\begin{equation}
    \eta_{yzx} = \frac{1}{2} \left[ D^{\alpha_2}_x [D^{\alpha_1}_y w_0] - D^{\alpha_2}_x \theta_y \right] = -\frac{1}{2} D^{\alpha_2}_x \theta_y
\end{equation}
\begin{equation}
    \eta_{yzy} = \frac{1}{2} \left[ D^{\alpha_2}_y [D^{\alpha_1}_y w_0] - D^{\alpha_2}_y \theta_y \right] = -\frac{1}{2} D^{\alpha_2}_y \theta_y
\end{equation}
\begin{equation}
    \eta_{yzz} = \frac{1}{2} \left[ D^{\alpha_2}_z [D^{\alpha_1}_y w_0] - D^{\alpha_2}_z \theta_y \right] = 0
\end{equation}
\end{subequations}
Assuming small displacement gradients $(O(\varepsilon))$, the strain-gradient terms except for those provided in Eq.~(\ref{eq: mindlin_exact_strain_grads}) are $O(\varepsilon^2)$. Thus, analogously to the arguments used in the development of the Timoshenko beam, for the normal strains the transverse strain-gradients are significant when compared to the in-plane gradients. When obtaining the solution via the weak form, the contribution of the in-plane strain-gradients of the normal strains can be ignored compared to the contribution of the non-zero strain-gradients in Eq.~(\ref{eq: mindlin_exact_strain_grads}). This observation can also be noted from results presented in \cite{jafari2016size}, where it is shown that ignoring the transverse strain-gradients of the normal strains leads to a significant change in the response of the structure, while the inclusion of the in-plane strain-gradients of the normal strains does not lead to a significant change in the response.

Using strains and strain-gradients in Eqs.~(\ref{eq: Mindlin_strains},\ref{eq: mindlin_exact_strain_grads}), the first variations of the potential energy, the kinetic energy and the work done by externally applied forces are obtained as:
\begin{subequations}
\label{eq: Mindlin_virtual_quantities}
\begin{equation}
\label{eq: mindlin_virtual_strain_energy}
\begin{split}
    \delta \mathcal{U} = \int_{\Omega_M} [ N_{xx}\delta D^{\alpha_1}_x u_0 + N_{yy}\delta D^{\alpha_1}_y v_0 + N_{xy}\delta \left(D^{\alpha_1}_y u_0 + D^{\alpha_1}_x v_0\right) + M_{xx}\delta D^{\alpha_1}_x \theta_x + M_{yy}\delta D^{\alpha_1}_y \theta_y  + \\
    ~M_{xy}\delta \left(D^{\alpha_1}_y \theta_x + D^{\alpha_1}_x \theta_y\right) +  Q_{xz} \delta(D^{\alpha_1}_x w_0-\theta_x) + Q_{yz} \delta(D^{\alpha_1}_y w_0-\theta_y) + \overline{N}_{xxz} \delta D^{\alpha_1}_x \theta_x + \\
   \overline{N}_{yyz}\delta D^{\alpha_1}_y \theta_y + \overline{N}_{xyz}\delta \left(D^{\alpha_1}_y \theta_x + D^{\alpha_1}_x \theta_y \right) + \overline{N}_{xzx}\delta D^{\alpha_2}_x \theta_x + \overline{N}_{xzy}\delta D^{\alpha_2}_y \theta_x + \\ 
   \overline{N}_{yzx}\delta D^{\alpha_2}_x \theta_y + \overline{N}_{yzy}\delta D^{\alpha_2}_y \theta_y] \mathrm{d} \Omega_M
\end{split}
\end{equation}
\begin{equation}
\label{eq: mindlin_virtual_work}
    \delta V = \int_{\Omega_M} \left[ F_x\delta u_0 + F_y \delta v_0 + F_z \delta w_0 + M_{\theta_x} \delta \theta_x + M_{\theta_y} \delta \theta_y \right]\mathrm{d} \Omega_M
\end{equation}
\begin{equation}
\label{eq: mindlin_virtual_kinetic_energy}
    \delta T = \int_{\Omega_M} \bigg\{\int_{-\frac{h}{2}}^{\frac{h}{2}} \rho \left[ \big(\dot{u}_0 - z\dot{\theta}_x \big)\big(\delta\dot{u}_0 - z\delta\dot{\theta}_x \big) + \big(\dot{v}_0 - z\dot{\theta}_y \big) \big(\delta\dot{v}_0 - z\delta\dot{\theta}_y \big) +\dot{w}_0 \delta \dot{w}_0 \right]\mathrm{d}z \bigg\} \mathrm{d} \Omega_M
\end{equation}
\end{subequations}
Note that $\mathrm{d}{\Omega_M} = \mathrm{d}x\mathrm{d}y$ for a rectangular plate. $\{F_x,F_y,F_z\}$ are the external loads applied in the $\hat{x}$, $\hat{y}$, and $\hat{z}$ directions, respectively. $\{M_{\theta_x},M_{\theta_y}\}$ are the external moments applied about the $\hat{y}$ and $\hat{x}$ axes, respectively. The different stress, moment, and higher-order stress resultants in Eq.~(\ref{eq: mindlin_virtual_kinetic_energy}) extend directly from Eq.~(\ref{eq: Timo_stress_resultants}).

The f-FEM for the Mindlin plates extends directly from the f-FEM formulation briefly reviewed in \S\ref{sec: Timoshenko_beams}. We also highlight that the f-FEM for fractional-order Mindlin plates can also be found in \cite{patnaik2020plates,patnaik2020geometrically}. Thus, for the sake of brevity, we do not provide provide all the details but we highlight the additional contributions following from the nonlocal strain-gradient terms. The expression for the stiffness matrix corresponding to the f-FEM for the Mindlin plate is:
\begin{equation}
\label{eq: Mindlin_FE_expression_virtual_strain}
    [K_M] = \int_{{\Omega_M}} [\tilde{B}_M(\bm{x})]^T [S_M] [\tilde{B}_M(\bm{x})] \mathrm{d}\Omega_M
\end{equation}
where $[S_M]$ is the constitutive matrix of the plate and the matrix $[\tilde{B}_M(\bm{x})]$ is given as:
\begin{equation}
\label{eq: Mindlin_FE_B_matrices_collection}
\begin{split}
        [\tilde{B}_M(\bm{x})] = \bigg[ \underbrace{[\tilde{B}^{\alpha_1}_{u_0,x}(\bm{x})]^T, [\tilde{B}^{\alpha_1}_{v_0,y}(\bm{x})]^T, \left[ [\tilde{B}^{\alpha_1}_{u_0,y}(\bm{x})] + [\tilde{B}^{\alpha_1}_{v_0,x}(\bm{x})] \right]^T, [\tilde{B}^{\alpha_1}_{\theta_x,x}(\bm{x})]^T, [\tilde{B}^{\alpha_1}_{\theta_y,y}(\bm{x})]^T}_{\text{Contributions from the nonlocal strains}} , \\ 
        \underbrace{\left[ [\tilde{B}^{\alpha_1}_{\theta_y,x}(\bm{x})] + [\tilde{B}^{\alpha_1}_{\theta_x,y}(\bm{x})] \right]^T, \left[\tilde{B}^{\alpha_1}_{w_0,x}(\bm{x})] - [{\mathbb{L}}^{(\theta_x)}(\bm{x})] \right]^T, \left[ [\tilde{B}^{\alpha_1}_{w_0,y}(\bm{x})] - [{\mathbb{L}}^{(\theta_y)}(\bm{x})] \right]^T}_{\text{Contributions from the nonlocal strains}}, \\
        \underbrace{[\tilde{B}^{\alpha_1}_{\theta_x,x}(\bm{x})]^T, [\tilde{B}^{\alpha_1}_{\theta_y,y}(\bm{x})]^T, \left[ [\tilde{B}^{\alpha_1}_{\theta_y,x}(\bm{x})] + [\tilde{B}^{\alpha_1}_{\theta_x,y}(\bm{x})] \right]^T}_{\text{Contributions from the nonlocal strain-gradients}}, 
        \\ 
        \underbrace{[\tilde{B}^{\alpha_2}_{\theta_x,x}(\bm{x})]^T, [\tilde{B}^{\alpha_2}_{\theta_x,y}(\bm{x})]^T, [\tilde{B}^{\alpha_2}_{\theta_x,x}(\bm{x})]^T, [\tilde{B}^{\alpha_2}_{\theta_x,y}(\bm{x})]^T}_{\text{Contributions from the nonlocal strain-gradients}} \bigg]
\end{split}
\end{equation}
The details of the fractional-order derivative matrices $[\tilde{B}^{\alpha_m}_{\square,r}(\bm{x})]$ can be found in SI and \cite{patnaik2020plates}.

%%%%%%%%%%%%%%%%%%%%%%%%%%%%%%%%%%%%%%%%%%%%%%%%%%%%%%%%%%%%%%%%%%%%%%%%%%%%%%%%%%%%%%%%%%%%%%%%%%%%%%
\subsection{Static response}
\label{ssec: Mindlin_static}
In this section, we analyze the static response of the Mindlin plate obtained via the fractional-order continuum formulation. For this purpose, the in-plane dimensions of the plate were fixed to be $L_M=1$m and $B_M=1$m and the thickness of the plate was taken to be $h_M=0.1$m $(=L/10)$. The simplified constitutive relations given in Eq.~(\ref{eq: simplfied_CE}) were used in this study. The material was assumed isotropic with an elastic modulus $E=30$GPa, Poisson's ratio $\nu=0.3$ and density $\rho=2700$kg/m${}^3$. Further, we have assumed a symmetric and isotropic horizon of nonlocality for points sufficiently inside the domain of the plate, that is $l_{A_\square} = l_{B_\square} = l_f, \square\in\{x,y\}$. For points located close to a boundary, the length scales were truncated as shown in Fig.~(\ref{fig: FCM}). 

We analyzed the static response of the plate subject to a UDTL of magnitude $F_z=10^{7}\text{Pa}$ for two different kinds of boundary conditions: the plate clamped at all the edges (CCCC) and the plate simply supported at all ts edges (SSSS) for different combinations of the fractional model parameters. For each boundary condition, we obtained the response of the plate for the following different cases:
\begin{itemize}
    \item Case 1: the fractional-orders $\alpha_1$ and $\alpha_2$ were varied within the range [0.5,1] for fixed values of the nonlocal horizon length $l_f=0.5\text{m}~(=L/2)$. For this case, the microstructural length was chosen as $l_*=0.02\text{m}~(=L/50)$. 
    
    \item Case 2: the nonlocal horizon $l_f$ was varied in $[0.5,1]\text{m}~ (=[L/2,L])$ and the microstructure length $l_*$ was varied in $[0.01,0.05]\text{m}~(=[L/100,L/20])$, for fixed values of the fractional-orders $\alpha_1=\alpha_2=0.8$ 
\end{itemize}

\begin{figure}[h!]
    \centering
    \includegraphics[width=\textwidth]{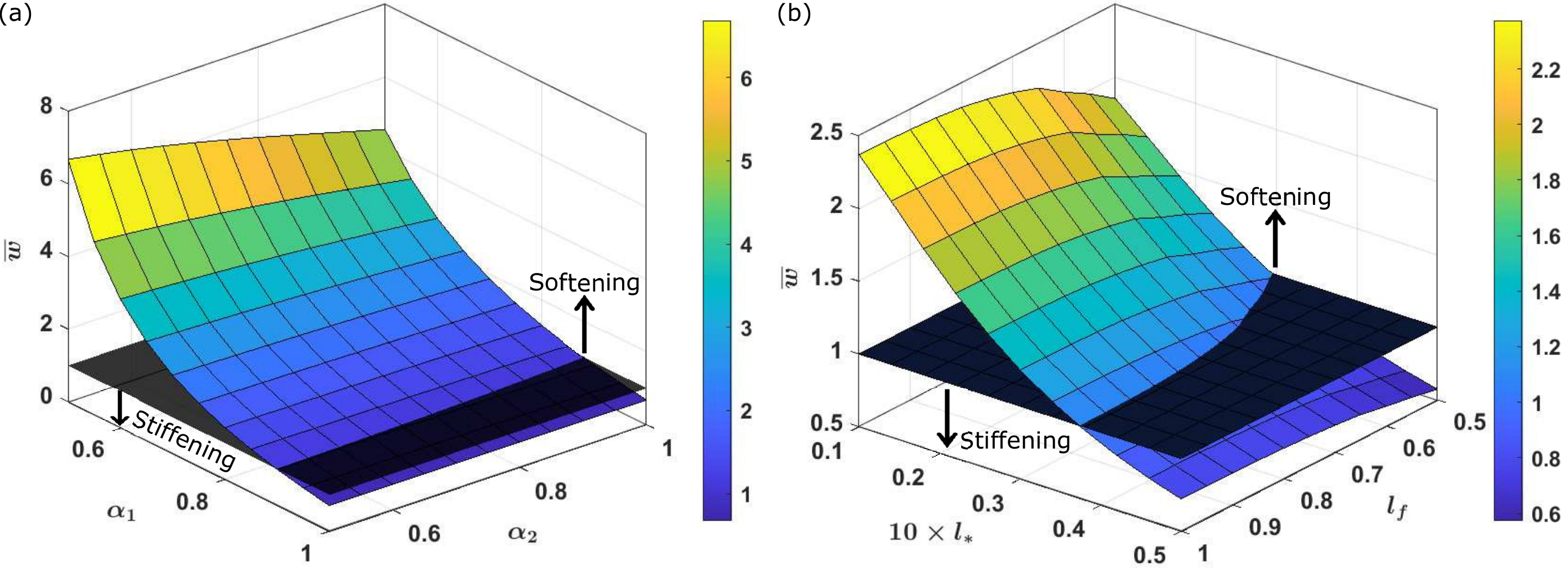}
    \caption{Non-dimensional transverse displacement at the center point of the Mindlin plate clamped at all its edges for different values of (a) the fractional orders and (b) the length scales.}
    \label{fig: Mindlin_CC}
\end{figure}

\begin{figure}[h!]
    \centering
    \includegraphics[width=\textwidth]{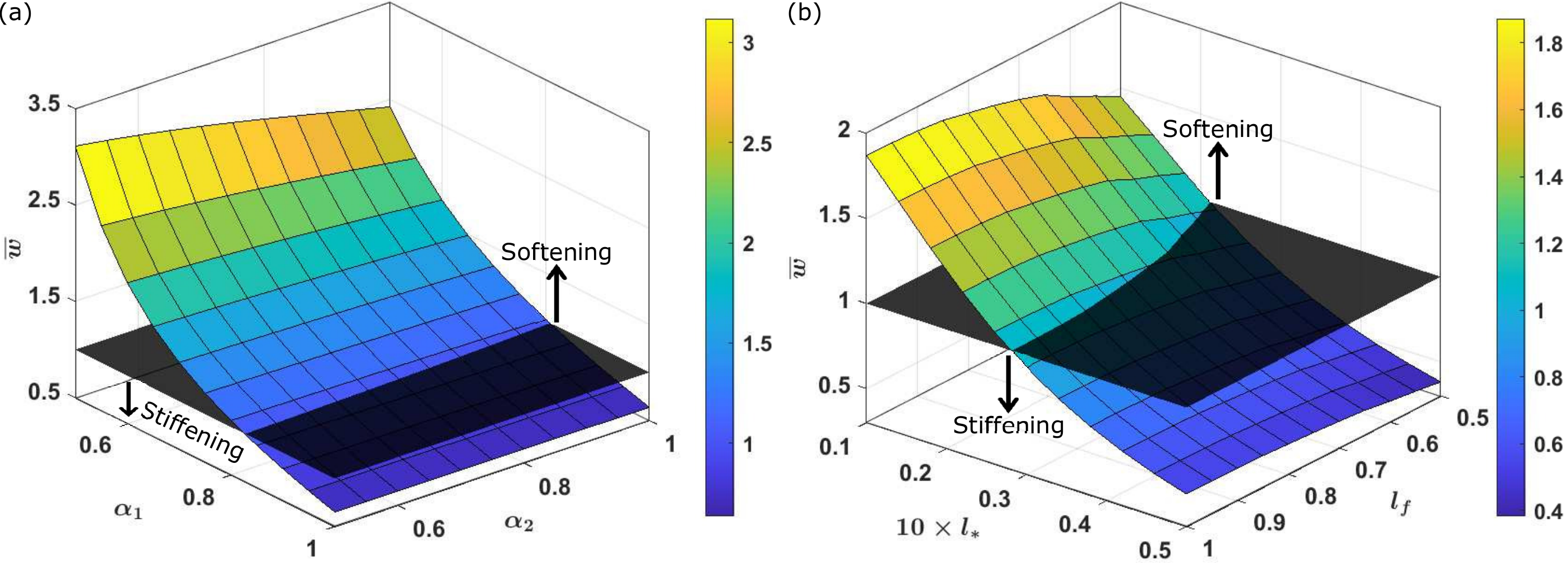}
    \caption{Non-dimensional transverse displacement at the center point of the Mindlin plate simply-supported at all its edges parameterized for different values of (a) the fractional orders and (b) the length scales.}
    \label{fig: Mindlin_SS}
\end{figure}

The numerical results, in terms of the maximum transverse displacement (obtained at the mid-point of the plate), are presented in Fig.~(\ref{fig: Mindlin_CC}) and Fig.~(\ref{fig: Mindlin_SS}) for the CCCC plate and the SSSS plate, respectively. Further, similar to the Timoshenko beam, the transverse displacement obtained for each case ($\overline{w}$) is non-dimensionalized against the maximum transverse displacement obtained for a classical Mindlin plate without nonlocality or strain-gradient effects. The maximum transverse displacement obtained for the classical CCCC plate was ${w}_0=0.55\times 10^{-2}\text{m}$ and for the classical SSSS plate was ${w}_0=1.55\times 10^{-2}\text{m}$. As evident from the Figs.~(\ref{fig: Mindlin_CC},\ref{fig: Mindlin_SS}), the fractional-order continuum formulation is able to model both stiffening and softening response of the Mindlin plate with respect to the classical formulation. The conclusions noted for the Timoshenko beam directly extend to the Mindlin plate. More specifically, the plate exhibits a stiffened response with increasing values of $\alpha_1$, $\alpha_2$ and $l_*$ and softened response with an increasing value of $l_f$ (see Fig.~(\ref{fig: parameter_map})).

%%%%%%%%%%%%%%%%%%%%%%%%%%%%%%%%%%%%%%%%%%%%%%%%%%%%%%%%%%%%%%%%%%%%%%%%%%%%%%%%%%%%%%%%%%%%%%%%%%%%%%
\subsection{Free vibration response}
\label{ssec: Mindlin_free_vibration}
In the following, we present the results capturing the effect of the different fractional model parameters on the natural frequency of transverse vibrations of the Mindlin plates. The material properties, loading conditions, boundary conditions and the range of the different fractional model parameters are the same as chosen for the static analysis of the Mindlin plate in \S\ref{ssec: Mindlin_static}. The results for the CCCC plate and the SSSS plate are presented in Figs.~(\ref{fig: Mindlin_CC},\ref{fig: Eigen_Mindlin_SS}), respectively, in terms of the non-dimensionalized natural frequency $\overline{\omega}_0$. Similar to the analysis in \ref{ssec: Timoshenko_free_vibration}, the non-dimensionalized natural frequency is obtained by dividing the natural frequency of the fractional-order plate with the natural frequency of the classical Mindlin plate for the specific boundary condition. The natural frequency obtained for the classical CCCC plate was ${\omega}_0=522\text{Hz}$ and for the classical SSSS plate was ${\omega}_0=306\text{Hz}$. As evident from the Figs.~(\ref{fig: Mindlin_CC},\ref{fig: Eigen_Mindlin_SS}), the conclusions presented in \S\ref{ssec: Timoshenko_static} on the specific effects of the different fractional model parameters, hold true for the free vibration response of the Mindlin plates.

\begin{figure}[h!]
    \centering
    \includegraphics[width=\textwidth]{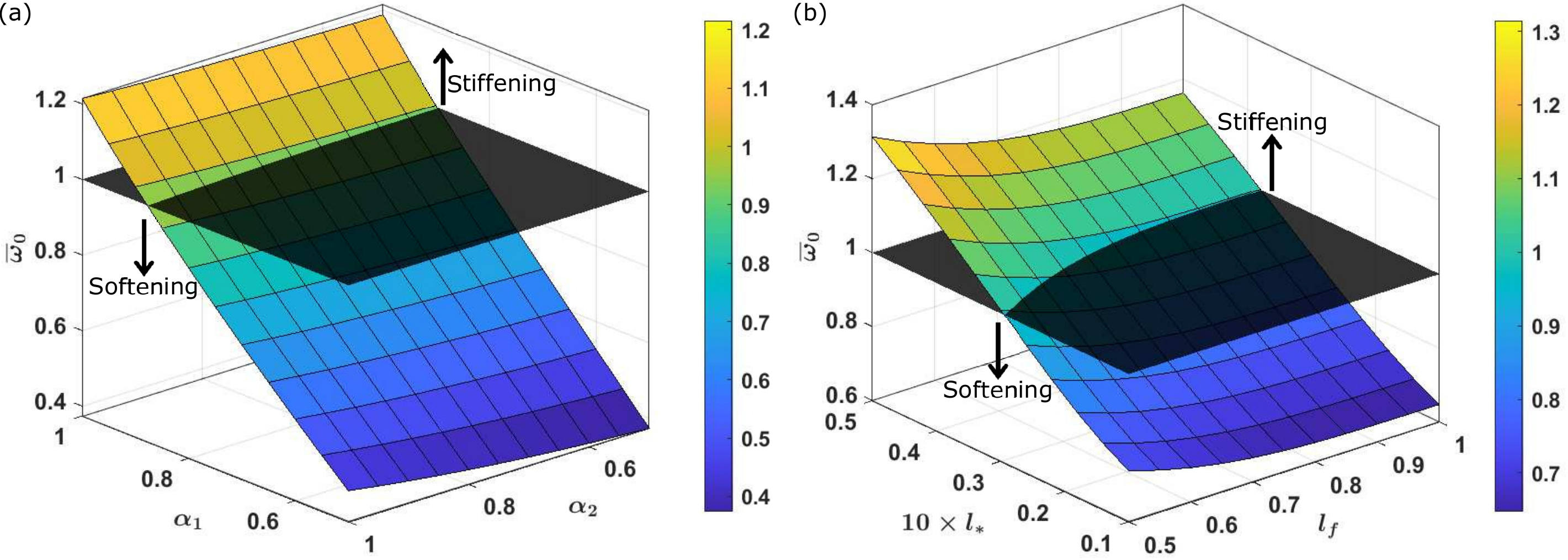}
    \caption{Non-dimensionalized natural frequency of the Mindlin plate clamped at its edges parameterized for different values of (a) the fractional orders and (b) the length scales.}
    \label{fig: Eigen_Mindlin_CC}
\end{figure}

\begin{figure}[h!]
    \centering
    \includegraphics[width=\textwidth]{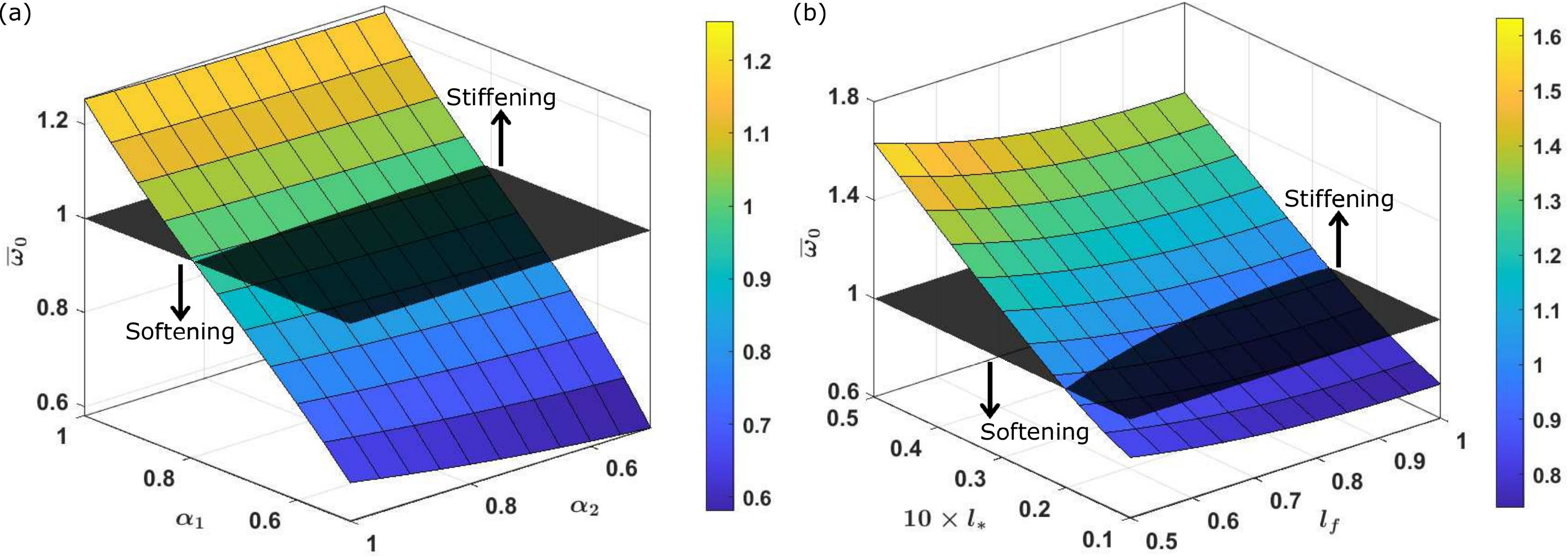}
    \caption{Non-dimensionalized natural frequency of the Mindlin plate simply supported at its boundaries parameterized for different values of (a) the fractional orders and (b) the length scales.}
    \label{fig: Eigen_Mindlin_SS}
\end{figure}

%%%%%%%%%%%%%%%%%%%%%%%%%%%%%%%%%%%%%%%%%%%%%%%%%%%%%%%%%%%%%%%%%%%%%%%%%%%%%%%%%%%%%%%%%%%%%%%%%%%%%%
\section{Conclusions}
The present study leveraged the fractional-order mechanics framework to develop a unified approach to nonlocal elasticity that combines the characteristics of both integral and gradient based classical formulations. More specifically, the differ-integral nature of fractional order operators was exploited to formulate a nonlocal continuum theory capable of modeling both stiffening and softening responses in structures exhibiting size-dependent effects. The fractional-order formulation was derived by the continualization of the Lagrangian of a 1D lattice subject to long-range cohesive interactions. Then, the governing equations corresponding to a 3D continuum were derived using variational principles. The resulting nonlocal theory is frame-invariant and causal. Contrary to classical integral formulations, the fractional-order formulation of a nonlocal continuum leads to positive definite systems with well-posed governing equations. Particularly remarkable is the ability of the fractional-order continuum model to capture anomalous attenuation and dispersion without having to incorporate inertia gradients in the governing equations; otherwise needed in classical strain-gradient formulations. Consequently, the fractional theory is well suited to capture nonlocality, scale effects, and medium heterogeneity in structural problems. The ability of the fractional-order formulation to model both stiffening and softening response was exemplified by performing both static and free vibration analysis of Timoshenko beams and Mindlin plates. In conclusion, the formulation and the results presented in the study illustrated several unique features of fractional calculus and suggested that this mathematical tool could play a critical role in the development of unified and comprehensive simulation tools for modeling the response of complex nonlocal structures.

\section{Acknowledgements}
The authors gratefully acknowledge the financial support of the Defense Advanced Research Project Agency (DARPA) under grant \#D19AP00052, and of the National Science Foundation (NSF) under grants MOMS \#1761423 and DCSD \#1825837. The content and information presented in this manuscript do not necessarily reflect the position or the policy of the government. The material is approved for public release; distribution is unlimited.

%\putbib[report]
%\end{bibunit}
\bibliographystyle{unsrt}
\bibliography{report}

\clearpage

{\noindent\underline{\LARGE{\textbf{Supplementary Information}}}}
\setcounter{equation}{0}
\setcounter{page}{1}
\setcounter{section}{0}
\setcounter{figure}{0}
\renewcommand{\theequation}{S\arabic{equation}}
\renewcommand{\thepage}{S\arabic{page}}
\renewcommand{\thesection}{S\arabic{section}}
\renewcommand{\thefigure}{S\arabic{figure}}
%\begin{bibunit}[unsrt]

\section{Derivation of the strong form of the 3D governing equations}
\label{sec: Appendix_A}
In the following, we have provided the detailed steps adopted in obtaining the first variation of the potential energy. The first variation of the potential energy follows from Eqs.~(21,25) as:
\begin{equation}
\label{eq: derivation_step1}
    \mathcal{U} = \int_{\Omega} \left[  {\sigma}_{ij} \delta {\epsilon}_{ij} + {\tau}_{ijk} \delta {\eta}_{ijk} \right] \mathrm{d}\mathbb{V}
\end{equation}
By using the kinematic relations in Eqs.~(22,24) and the symmetry of the stress and higher-order stress tensor, Eq.~(\ref{eq: derivation_step1}) is expressed as:
\begin{equation}
\label{eq: derivation_step2}
    \mathcal{U} = \int_{\Omega} \big[  \underbrace{{\sigma}_{ij} \delta (D^{\alpha_1}_{x_j} u_i)}_{\mathcal{I}_1} + \underbrace{\tau_{ijk} \delta \big[ D^{\alpha_2}_{x_k} (D^{\alpha_1}_{x_j} u_i)}_{\mathcal{I}_2} \big] \big] \mathrm{d}\mathbb{V}
\end{equation}
where the displacement field is assumed as: $\bm{u} = u \hat{x} + v \hat{y} + w \hat{z}$. The domain $\Omega=[0,L]\times[0,B]\times[0,H]$ is illustrated in Fig.~(3) of the manuscript.\\

\noindent\textbf{Simplification of $\mathcal{I}_1$:} The term $\mathcal{I}_1$ in Eq.~(\ref{eq: derivation_step2}) can be expanded as:
\begin{equation}
\label{eq: stress_step1}
\begin{split}
    \mathcal{I}_{1} = \int_{\Omega} \big[  \underbrace{{\sigma}_{xx} \delta (D^{\alpha_1}_{x} u)}_{\mathcal{I}_{11}} + {\sigma}_{yy} \delta (D^{\alpha_1}_{y} v) + & {\sigma}_{zz} \delta (D^{\alpha_1}_{z} w) + {\sigma}_{xy} \delta (D^{\alpha_1}_{y} u + D^{\alpha_1}_{x} v) + \\ & {\sigma}_{xz} \delta (D^{\alpha_1}_{z} u + D^{\alpha_1}_{x} w) + {\sigma}_{yz} \delta (D^{\alpha_1}_{z} v + D^{\alpha_1}_{y} w) \big] \mathrm{d}\mathbb{V}
\end{split}
\end{equation}
From the definition of the RC derivative in Eq.~(26), the first term within the above integral is expressed as:
\begin{equation}
\label{eq: stress_step2}
\begin{split}
    \mathcal{I}_{11} = \int_{0}^{H} \int_{0}^{B}\left[ \frac{1}{2}\Gamma(2-\alpha_1) \left\{ l_{A_x}^{\alpha_1-1} \int_0^L \sigma_{xx} \left[ {}^{~~~~~C}_{x-l_{A_x}}D^{\alpha_1}_{x}(\delta u) \right] \mathrm{d}x -
    l_{B_x}^{\alpha-1} \int_0^L \sigma_{xx}\left[{}^{C}_{x}D^{\alpha_1}_{x+l_{B_x}}(\delta u)\right]\mathrm{d}x \right\} \right] \mathrm{d}y \mathrm{d} z
\end{split}
\end{equation}
Using the definitions of the left- and right-handed Caputo derivatives we obtain:
\begin{subequations}
\begin{equation}
\label{eq: left_caputo_axial}
     \int_0^L \sigma_{xx} \left[ {}^{~~~~~C}_{x-l_{A_x}}D^{\alpha_1}_{x}(\delta u)\right]\mathrm{d}x = \frac{1}{\Gamma(1-\alpha_1)}\int_0^L \sigma_{xx} \left[\int_{x - l_{A_x}}^{x}\left(x-x^\prime\right)^{-\alpha_1} [D^1_{x^\prime} (\delta u)] \mathrm{d}x^\prime\right]\mathrm{d}x
\end{equation}
\begin{equation}
\label{eq: right_caputo_axial}
    \int_0^L \sigma_{xx}\left[{}^{C}_{x}D^{\alpha_1}_{x+l_{B_x}}(\delta u)\right]\mathrm{d}x = -\frac{1}{\Gamma(1-\alpha_1)}\int_0^L \sigma_{xx} \left[ \int_{x}^{x+l_{B_x}} \left(x^\prime - x \right)^{-\alpha_1} [D^1_{x^\prime} (\delta u)] \mathrm{d}x^\prime \right] \mathrm{d}x
\end{equation}
\end{subequations}
where $x^\prime$ is a dummy spatial variable along the $\hat{x}$ direction used for the convolution. The above integrals are further evaluated using integration by parts in order to transfer the derivative from independent variable (displacement field) to the secondary variable (stress). This leads to the following:
\begin{subequations}
\begin{equation}
    \int_0^L \sigma_{xx} \left[ {}^{~~~~~C}_{x-l_{A_x}}D^{\alpha_1}_{x}(\delta u)\right]\mathrm{d}x = \frac{1}{\Gamma(1-\alpha_1)} \int_0^L \frac{\mathrm{d}\delta u(x^\prime)}{\mathrm{d}x^\prime} \left[\int_{x^\prime}^{x^\prime+l_A} \left(x-x^\prime\right)^{-\alpha_1} \sigma_{xx} \mathrm{d}x \right] \mathrm{d}x^\prime
\end{equation}
\begin{equation}
   \int_0^L \sigma_{xx}\left[{}^{C}_{x}D^{\alpha_1}_{x+l_{B_x}}(\delta u)\right]\mathrm{d}x = -\frac{1}{\Gamma(1-\alpha_1)} \int_0^L \frac{\mathrm{d}\delta u(x^\prime)}{\mathrm{d}x^\prime} \left[\int_{x^\prime-l_{B_x}}^{x^\prime}\left(x^\prime-x\right)^{-\alpha_1} \sigma_{xx} \mathrm{d}x \right] \mathrm{d}x^\prime
\end{equation}
\end{subequations}
Using the definitions for left- and right- fractional integrals in the above results, we obtain:
\begin{subequations}
\begin{equation}
\label{eq: left_caputo_axial_intpart1}
    \int_0^L \sigma_{xx} \left[ {}^{~~~~~C}_{x-l_{A_x}}D^{\alpha_1}_{x}(\delta u)\right]\mathrm{d}x =\int_0^L [D^1_{x^\prime} (\delta u)] \left( {}_{x}I^{1-\alpha_1}_{x+l_{A_x}} [\sigma_{xx}] \right) \mathrm{d}x
\end{equation}
\begin{equation}
\label{eq: right_caputo_axial_intpart1}
   \int_0^L \sigma_{xx}\left[{}^{C}_{x}D^{\alpha_1}_{x+l_{B_x}}(\delta u)\right]\mathrm{d}x = -\int_0^L [D^1_{x^\prime} (\delta u)] \left( {}_{x-l_{B_x}}I^{1-\alpha_1}_{x} [\sigma_{xx}] \right) \mathrm{d}x
\end{equation}
\end{subequations}
Repeating the integration by parts and substituting the resulting expressions in Eq.~(\ref{eq: stress_step2}) we obtain:
\begin{equation}
\label{eq: stress_step3}
\begin{split}
    \mathcal{I}_{11} = \int_{0}^{H} \int_{0}^{B} \left[ \frac{1}{2}\Gamma(2-\alpha_1) \left\{ l_{A_x}^{\alpha_1-1} \left. \left[ \left( {}_{x}I^{1-\alpha_1}_{x + l_{A_x}} \sigma_{xx} \right) \delta u \right]\right\vert_{x=0}^{x=L}-
    l_{B_x}^{\alpha_1-1} \left. \left[ \left( {}_{x-l_{B_x}}I^{1-\alpha_1}_{x} \sigma_{xx} \right) \delta u \right] \right\vert_{x=0}^{x=L} \right\} \right] \mathrm{d}y \mathrm{d} z + \\
    \int_{0}^{H} \int_{0}^{B} \left[ \frac{1}{2}\Gamma(2-\alpha_1) \left\{l_{A_x}^{\alpha_1-1} \int_0^L \left[ {}^{RL}_{x} D^{\alpha_1}_{x + l_{A_x}} \sigma_{xx} \right] \delta u \mathrm{d}x - l_{B_x}^{\alpha_1-1} \int_0^L \left[{}^{~~~~RL}_{x-l_{B_x}}D^{\alpha_1}_{x} \sigma_{xx} \right]\delta u \mathrm{d}x \right\} \right] \mathrm{d}y \mathrm{d} z 
\end{split}
\end{equation}
Now by using the definitions for the Riesz fractional integral in Eq.~(34) and the Riesz Riemann Liouville derivative in Eq.~(36), the above expression can be simplified as:
\begin{equation}
\label{eq: stress_step4}
\begin{split}
    \mathcal{I}_{11} = - \underbrace{\int_{0}^{H} \int_{0}^{B} \int_{0}^{L} \left[ \mathfrak{D}^{\alpha_1}_x \sigma_{xx} \right] \delta u \mathrm{d}x \mathrm{d}y \mathrm{d} z}_{\text{Indicates integral over the volume $\Omega$}} + \underbrace{\int_{0}^{H} \int_{0}^{B}  \left. \left[ \left( I^{1-\alpha_1}_x \sigma_{xx} \right) \delta u \right] \right\vert_{x=0}^{x=L}  \mathrm{d}y \mathrm{d} z}_{\text{Indicates integration over surfaces $\partial\Omega_x$}} \\ = 
    - \int_{\Omega} \left[ \mathfrak{D}^{\alpha_1}_x \sigma_{xx} \right] \delta u \mathrm{d}\mathbb{V} + \int_{\partial \Omega_x} \left[ I^{1-\alpha_1}_x \sigma_{xx} \right] \delta u \mathrm{d}\mathbb{A}
    \end{split}
\end{equation}
Similar variational simplifications have also been carried out for 1D and 2D BVPs in \cite{patnaik2019FEM,patnaik2020plates}. Using the above outlined steps, the remaining terms in $\mathcal{I}_1$ can be simplified as:
\begin{subequations}
\label{eq: stress_terms_simplified}
\begin{equation}
\begin{split}
    \int_{\Omega} \big[  {\sigma}_{xx} \delta (D^{\alpha_1}_{x} u) + {\sigma}_{xy} \delta (D^{\alpha_1}_{x} v) + {\sigma}_{xz} \delta (D^{\alpha_1}_{x} w) \big] \mathrm{d}\mathbb{V} =
    -\int_{\Omega} \left[ \left(\mathfrak{D}^{\alpha_1}_x \sigma_{xx} \right) \delta u + \left(\mathfrak{D}^{\alpha_1}_x \sigma_{xy} \right) \delta v + \left(\mathfrak{D}^{\alpha_1}_x \sigma_{xz} \right) \delta w \right] \mathrm{d}\mathbb{V} + \\
    \int_{\partial \Omega_x} \left[ \left( I^{1-\alpha_1}_x \sigma_{xx} \right) \delta u + \left( I^{1-\alpha_1}_x \sigma_{xy} \right) \delta v + \left( I^{1-\alpha_1}_x \sigma_{xz} \right) \delta w  \right] \mathrm{d}\mathbb{A} 
\end{split}
\end{equation}
\begin{equation}
\begin{split}
    \int_{\Omega} \big[  {\sigma}_{xy} \delta (D^{\alpha_1}_{y} u) + \sigma_{yy} \delta (D^{\alpha_1}_{y} v) + {\sigma}_{yz} \delta (D^{\alpha_1}_{y} w) \big] \mathrm{d}\mathbb{V} =
    -\int_{\Omega} \left[ \left(\mathfrak{D}^{\alpha_1}_y \sigma_{xy} \right) \delta u + \left(\mathfrak{D}^{\alpha_1}_y \sigma_{yy} \right) \delta v + \left(\mathfrak{D}^{\alpha_1}_y \sigma_{yz} \right) \delta w \right] \mathrm{d}\mathbb{V} + \\
    \int_{\partial \Omega_y} \left[ \left( I^{1-\alpha_1}_y \sigma_{xy} \right) \delta u + \left( I^{1-\alpha_1}_y \sigma_{yy} \right) \delta v + \left( I^{1-\alpha_1}_y \sigma_{yz} \right) \delta w  \right] \mathrm{d}\mathbb{A} 
\end{split}
\end{equation}
\begin{equation}
\begin{split}
    \int_{\Omega} \big[  {\sigma}_{xz} \delta (D^{\alpha_1}_{z} u) + {\sigma}_{yz} \delta (D^{\alpha_1}_{z} v) + {\sigma}_{zz} \delta (D^{\alpha_1}_{z} w) \big] \mathrm{d}\mathbb{V} =
    -\int_{\Omega} \left[ \left(\mathfrak{D}^{\alpha_1}_z \sigma_{xz} \right) \delta u + \left(\mathfrak{D}^{\alpha_1}_z \sigma_{yz} \right) \delta v + \left(\mathfrak{D}^{\alpha_1}_z \sigma_{zz} \right) \delta w \right] \mathrm{d}\mathbb{V} + \\
    \int_{\partial \Omega_z} \left[ \left( I^{1-\alpha_1}_z \sigma_{xz} \right) \delta u + \left( I^{1-\alpha_1}_x \sigma_{yz} \right) \delta v + \left( I^{1-\alpha_1}_z \sigma_{zz} \right) \delta w  \right] \mathrm{d}\mathbb{A} 
\end{split}
\end{equation}
\end{subequations}
Combining all the terms in the above equation and using the definition of the integral operator $\bm{I}^{1-\alpha_1}_{\bm{\hat{n}}}(\cdot)$ in Eq.~(33) and the definition of the gradient operator $\bm{\widetilde{\nabla}}^{\alpha}(\cdot)$ in Eq.~(35), it immediately follows that:
\begin{equation}
    \label{eq: final_stress_simplified}
    \mathcal{I}_1 = -\int_{\Omega} \left[ \left(  \bm{\widetilde{\nabla}}^{\alpha} \cdot \bm{\sigma} \right) \cdot \delta \bm{u} \right] \mathrm{d}\mathbb{V} + \\
    \int_{\partial \Omega} \left[ \left(  \bm{I}^{1-\alpha_1}_{\bm{\hat{n}}} \cdot \bm{\sigma} \right) \cdot \delta \bm{u} \right] \mathrm{d}\mathbb{A} 
\end{equation}
It remains to simplify the variation of the strain-gradient energy contributions in the potential energy, i.e., the term $\mathcal{I}_2$ in Eq.~(\ref{eq: derivation_step2}).\\

\noindent\textbf{Simplification of $\mathcal{I}_2$:} The term $\mathcal{I}_2$ in Eq.~(\ref{eq: derivation_step2}) is given as:
\begin{equation}
\label{eq: strain_grad_step1}
    \mathcal{I}_2 = \int_{\Omega} \tau_{ijk} \delta \big[ D^{\alpha_2}_{x_k} (D^{\alpha_1}_{x_j} u_i) \big] \mathrm{d}\mathbb{V}
\end{equation}
To simplifify the expression above, we consider two cases: (C1) $j=k$ and (C2) $j\neq k$.\\

\noindent\textbf{Case C1:} For the case C1 when $j=k$, following the steps through Eqs.~(\ref{eq: stress_step2}-\ref{eq: stress_step4}), it follows immediately that:
\begin{equation}
\label{eq: strain_grad_case1_step1}
    \int_{\Omega} \tau_{ikk} \delta \big[ D^{\alpha_2}_{x_k} (D^{\alpha_1}_{x_k} u_i) \big] \mathrm{d}\mathbb{V} = 
    - \int_{\Omega} \underbrace{\left[ \mathfrak{D}^{\alpha_2}_{x_k} \tau_{ikk} \right] \big[ \delta (D^{\alpha_1}_{x_k} u_i) \big]}_{\mathcal{I}_{21}} \mathrm{d}\mathbb{V} + \int_{\partial \Omega_{x_k}} \underbrace{\left[ I^{1-\alpha_2}_{x_k} \tau_{ikk} \right] \big[ \delta (D^{\alpha_1}_{x_k} u_i) \big]}_{\mathcal{I}_{22}} \mathrm{d}\mathbb{A}
\end{equation}
In the above equation, the normal gradient indicated within the term ${\mathcal{I}_{22}}$ varies independently of the variation of $\bm{u}$ (since $j=k$), analogously to classical strain-gradient formulations.  Hence, the term ${\mathcal{I}_{22}}$ does not need to be simplified any further. The term ${\mathcal{I}_{21}}$ is simplified by retracing the steps through Eqs.~(\ref{eq: stress_step2}-\ref{eq: stress_step4}). Consequently, the above equation is further simplified as:
\begin{equation}
\label{eq: strain_grad_case1_step2}
\begin{split}
   \int_{\Omega} \tau_{ikk} \delta \big[ D^{\alpha_2}_{x_k} (D^{\alpha_1}_{x_k} u_i) \big] \mathrm{d}\mathbb{V} = 
    \int_{\Omega} \left[ \mathfrak{D}^{\alpha_1}_{x_k} \left( \mathfrak{D}^{\alpha_2}_{x_k} \tau_{ikk} \right) \right] \delta u_i \mathrm{d}\mathbb{V} - \int_{\partial \Omega_{x_k}} \left[ I^{1-\alpha_1}_{x_k} \left( \mathfrak{D}^{\alpha_2}_{x_k} \tau_{ikk} \right) \right] \delta u_i \mathrm{d}\mathbb{A} + \\ \int_{\partial \Omega_{x_k}} \left[ I^{1-\alpha_2}_{x_k} \tau_{ikk} \right] \big[ \delta (D^{\alpha_1}_{x_k} u_i) \big]\mathrm{d}\mathbb{A} 
\end{split}
\end{equation}

\noindent\textbf{Case C2:} following the steps outlined by Eqs.~(\ref{eq: stress_step2}-\ref{eq: stress_step4}), it follows that:
\begin{equation}
\label{eq: strain_grad_case2_step1}
    \int_{\Omega} \tau_{ijk} \delta \big[ D^{\alpha_2}_{x_k} (D^{\alpha_1}_{x_j} u_i) \big] \mathrm{d}\mathbb{V} = 
    - \int_{\Omega} {\left[ \mathfrak{D}^{\alpha_2}_{x_k} \tau_{ijk} \right] \big[ \delta (D^{\alpha_1}_{x_j} u_i) \big]} + \int_{\partial \Omega_{x_k}} {\left[ I^{1-\alpha_2}_{x_k} \tau_{ijk} \right] \big[ \delta (D^{\alpha_1}_{x_j} u_i) \big]} \mathrm{d}\mathbb{A}
\end{equation}
Note that unlike Case 1, the surface integral has to be further evaluated to relieve the variation of the displacement $u_i$ of the fractional-order gradient. The  above expression is simplified by again retracing the steps through Eqs.~(\ref{eq: stress_step2}-\ref{eq: stress_step4}), to obtain:
\begin{equation}
\label{eq: strain_grad_case2_step2}
\begin{split}
   \int_{\Omega} \tau_{ijk} \delta \big[ D^{\alpha_2}_{x_k} (D^{\alpha_1}_{x_j} u_i) \big] \mathrm{d}\mathbb{V} = 
    \int_{\Omega} \left[ \mathfrak{D}^{\alpha_1}_{x_j} \left( \mathfrak{D}^{\alpha_2}_{x_k} \tau_{ijk} \right) \right] \delta u_i \mathrm{d}\mathbb{V} - \int_{\partial \Omega_{x_j}} \left[ I^{1-\alpha_1}_{x_j} \left( \mathfrak{D}^{\alpha_2}_{x_k} \tau_{ijk} \right) \right] \delta u_i \mathrm{d}\mathbb{A} - \\ \int_{\partial \Omega_{x_k}} \left[ \mathfrak{D}^{\alpha_1}_{x_j} \left( I^{1-\alpha_2}_{x_k} \tau_{ijk} \right) \right] \delta u_i \mathrm{d}\mathbb{A} + \oint_{\overline{\Gamma}_{x_k}} \left[ I^{1-\alpha_1}_{x_j} \left( I^{1-\alpha_2}_{x_k} \tau_{ijk} \right) \right] \delta u_i\mathrm{d}\mathbb{A} 
\end{split}
\end{equation}\\

\noindent Combining all the terms obtained in Eqs.~(\ref{eq: strain_grad_case1_step2},\ref{eq: strain_grad_case2_step2}) the variation of the strain-gradient contributions to potential energy is obtained as:
\begin{equation}
    \label{eq: strain_grad_step2}
    \begin{split}
    \mathcal{I}_2 = 
    \int_{\Omega} \underbrace{\left[ \mathfrak{D}^{\alpha_1}_{x_j} \left( \mathfrak{D}^{\alpha_2}_{x_k} \tau_{ijk} \right) \right] \delta u_i}_{\mathcal{I}_2^{(1)}} \mathrm{d}\mathbb{V} - \int_{\partial \Omega_{x_j}} \underbrace{\left[ I^{1-\alpha_1}_{x_j} \left( \mathfrak{D}^{\alpha_2}_{x_k} \tau_{ijk} \right) \right] \delta u_i}_{\mathcal{I}_2^{(2)}} \mathrm{d}\mathbb{A} + \int_{\partial \Omega_{x_k}} \underbrace{\delta_{ij}\left[ I^{1-\alpha_2}_{x_k} \tau_{ijk} \right] \big[ \delta (D^{\alpha_1}_{x_j} u_i) \big]}_{\mathcal{I}^{(3)}_2} \mathrm{d}\mathbb{A} - \\  \int_{\partial \Omega_{x_k}} \underbrace{(1- \delta_{ij}) \left[ \mathfrak{D}^{\alpha_1}_{x_j} \left( I^{1-\alpha_2}_{x_k} \tau_{ijk} \right) \right] \delta u_i}_{\mathcal{I}_2^{(4)}} \mathrm{d}\mathbb{A} + \oint_{\overline{\Gamma}_{x_k}} \underbrace{(1 - \delta_{ij}) \left[ I^{1-\alpha_1}_{x_j} \left( I^{1-\alpha_2}_{x_k} \tau_{ijk} \right) \right] \delta u_i}_{\mathcal{I}_2^{(5)}} \mathrm{d}\mathrm{l}
    \end{split}
\end{equation}
where $\delta_{ij}$ represents the Kronecker delta function. The above equation can be represented as:
\begin{equation}
\label{eq: strain_grad_step3}
\begin{split}
\mathcal{I}_2 = \int_{\Omega} \underbrace{{\bm{\widetilde{\nabla}}}^{\alpha_1} \cdot \big( {\bm{\widetilde{\nabla}}}^{\alpha_2} \cdot \bm{\tau} \big) \cdot \delta \bm{u}}_{\mathcal{I}_2^{(1)}} \mathrm{d}\mathbb{V} 
- \int_{\partial \Omega} \underbrace{\left[ \bm{I}^{1-\alpha_1}_{\bm{\hat{n}}} \cdot \big( {\bm{\widetilde{\nabla}}}^{\alpha_2} \cdot \bm{\tau} \big) +  \bm{R} : \left( \bm{I}^{1-\alpha_2}_{\bm{\hat{n}}} \cdot \bm{\tau} \right) \otimes {\bm{\widetilde{\nabla}}}^{\alpha_1} \right] \cdot \delta \bm{u}}_{ \mathcal{I}_2^{(2)} + \mathcal{I}_2^{(4)} } \mathrm{d}\mathbb{A}  + \\
+\int_{\partial \Omega} \underbrace{\left[ \bm{I}^{1-\alpha_2}_{\bm{\hat{n}}} \otimes \bm{\hat{n}} \right] : \bm{\tau} \cdot \left[ \bm{\hat{n}} \cdot \left( {\delta} \bm{u} \otimes {\bm{{\nabla}}}^{\alpha_1} \right) \right]}_{\mathcal{I}_2^{(3)}} \mathrm{d}\mathbb{A}
+ \oint_{\overline{\Gamma}} \underbrace{ \left[\left[ \bm{I}^{1-\alpha_1}_{\bm{\hat{m}}} \cdot ( \bm{I}^{1-\alpha_2}_{\bm{\hat{n}}} \cdot \bm{\tau} ) \right] \right] \cdot \delta \bm{u}}_{\mathcal{I}_2^{(5)}} \mathrm{d}\mathrm{l}
\end{split}
\end{equation}
where we have indicated the correspondence between the indicial and vector notations. Recall from \S{3.2}, that the tensor $\bm{R}$ is the projector onto the surface $\partial\Omega$, $\bm{\hat{m}}$ is the co-normal vector at the edges and [[$\cdot$]] denotes difference of its argument across both sides of the edge $\overline{\Gamma}$. By combining Eqs.~(\ref{eq: final_stress_simplified}, \ref{eq: strain_grad_step3}), we obtain the first variation of the potential energy as given in Eq.~(31).

\section{Fractional-order finite element formulation}
In this section, we provide the key highlights of the f-FEM used to numerically simulate the fractional-order system. The details of the f-FEM are extensive and the interested reader can find these details in \cite{patnaik2019FEM,patnaik2020plates}. The f-FEM for the Timoshenko beam is formulated by obtaining a discretized form of the first variation of the Lagrangian of the beam. We start by deriving the discretized form of the potential energy. The stress and moment resultants in Eq.~(47) can be expressed as:
\begin{equation}
\label{eq: FE_Constitutive_expression}
    \left\{ N_{xx}, Q_{xz}, M_{xx}, \overline{N}_{xxz}, \overline{N}_{xzx} \right\}^T = [S_T]\left\{ D^{\alpha_1}_x u_0, D^{\alpha_1}_x w_0 - \theta_x, D^{\alpha_1}_x \theta_0, D^{\alpha_1}_x \theta_0,  D^{\alpha_2}_x \theta_0  \right\}^T
\end{equation}
where $[S_T]$ denotes the constitutive matrix of the beam. It is immediate that the approximation of the potential energy requires the approximation of the different fractional-order derivatives in Eq.~(\ref{eq: FE_Constitutive_expression}). 

For this purpose, the beam domain $\Omega_T=[0,L]$ is uniformly discretized into disjoint three-noded line elements. The vector containing the nodal degrees of freedom of the element $\Omega^e_i \in \Omega_T$ is denoted as $\{U^e_i\}$ while, the global degrees of freedom vector is denoted as $\{U\}$. The unknown displacement field variables at any point ${x} \in \Omega^e_i$ are evaluated by interpolating the corresponding nodal degrees of freedom of $\Omega^e_i$. For example, the axial displacement $u_0$ at a point $x\in\Omega^e_i$ can be obtained as:
\begin{equation}
    \label{eq: u_interpolation}
    u_0(x) = \left\{ \mathbb{L}^{(u_0)}_i (x) \right\} \{U^e_i\} 
\end{equation}
where, $\left\{ \mathbb{L}^{(u_0)}_i (x) \right\}$ contains the Lagrangian shape functions for three-noded 1D elements. The superscript in the row vector $\left\{ \mathbb{L}^{(u_0)}_i (x) \right\}$ indicates the specific displacement variable being interpolated which is $u_0$ in Eq.~(\ref{eq: u_interpolation}) and the subscript denotes the element number. The fractional-order derivative $D^\alpha_x u_0$ at the point $x$ is obtained as:
\begin{equation}
\label{eq: frac_der_u_final}
    D^{\alpha}_{x}\left[u_0(x)\right] = \left[\int_{x - l_{A_x}}^{x + l_{B_x}} \mathcal{K}(x, x^\prime, l_{A_x}, l_{B_x}, \alpha) [B_{u_0,x}(x')] [\tilde{\mathcal{C}}(x,x')] ~\mathrm{d}x'\right]\{U\} = [\tilde{B}^{\alpha}_{u_0,x}(x)]\{U\}
\end{equation}
where $x^\prime$ is a dummy variable used for convolution along the $\hat{x}$ axis. Note that $x'$ lies in the domain $(x - l_{A_x} , x + l_{B_x})$, which is the horizon of nonlocality at $x$. The remaining terms introduced in Eq.~(\ref{eq: frac_der_u_final}) are explained in the following. $\mathcal{K}(x, x^\prime, l_{A_x}, l_{B_x}, \alpha) [B_{u,x}(x')]$ denotes the kernel of the fractional-order derivative:
\begin{equation}
    \label{eq: frac_kernel_gen}
    \mathcal{K}(x, x^\prime, l_{A_x}, l_{B_x}, \alpha) = \begin{cases}
    \frac{1}{2} (1-\alpha) ~{l_{A_x}}^{\alpha-1} ~{|x - x^\prime|^{-\alpha}} & ~~ x^\prime \in (x-l_{A_x}, x)\\
    \frac{1}{2} (1-\alpha) ~{l_{B_x}}^{\alpha-1} ~{|x - x^\prime|^{-\alpha}} & ~~ x^\prime \in (x, x+ l_{B_x})
\end{cases}
\end{equation}
Note that the definition of $D^{\alpha}_{x}\left[u_0(x)\right]$ contains the integer-order derivative $D^1_{x^\prime}[u_0(x^\prime)]$. $D^1_{x^\prime}[u_0(x^\prime)]$ is evaluated at $x'$ in terms of the nodal displacement variables corresponding to the element $\Omega_p^e$, such that $x'\in \Omega_p^e$. Using Eq.~(\ref{eq: u_interpolation}), the integer-order derivative can be expressed as:
\begin{equation}
\label{eq: deriv_disp_rel}
    D^1_{x'} [u_0(x')] = [B_{u_0,x}(x')]\{U^e_p\} = \frac{\partial}{\partial x}  \left[ \left\{ \mathbb{L}^{(u_0)}_p ({x}') \right\} \right] \{U^e_p\}
\end{equation}
Further, $[\tilde{\mathcal{C}}(x,x')]$ is a connectivity matrix that is used to attribute the nonlocal contributions from the different elements in the horizon of $x$ to the corresponding nodes of those elements. In order to correctly account for these nonlocal contributions from the elements in the horizon, we transform the nodal values $\{U^e_{x'}\}$ into \{U\} using connectivity matrices in the following manner:
\begin{equation}
\label{eq: conversion_to_global_form}
    \{U^e_{x'}\} = [\tilde{\mathcal{C}}(x,x')]\{U\}
\end{equation}
The connectivity matrix $[\tilde{\mathcal{C}}(x,x')]$ is designed such that it is non-zero only if the point $x'$ lies in the nonlocal horizon of $x$. It is immediate to see that these matrices activate the contribution of the nodes enclosing $x'$ for the numerical evaluation of the convolution integral in Eq.~\eqref{eq: frac_der_u_final}. 

Following the above outlined procedure, the remaining fractional derivatives in Eq.~(\ref{eq: FE_Constitutive_expression}) are obtained as:
\begin{subequations}
\label{eq: FE_frac_derivatives_final}
\begin{equation}
    D^{\alpha}_{x}\left[\theta_0(x)\right] = [\tilde{B}^{\alpha}_{\theta_0,x}(x)] \{U\}
\end{equation}
\begin{equation}
   D^{\alpha}_{x}\left[w_0(x)\right] - \theta_x = \left[ [\tilde{B}^{\alpha}_{w_0,x}(x)] - [{\mathbb{L}}^{(\theta_0)}(x)] \right] \{U\}
\end{equation}
\end{subequations}
where $[{\mathbb{L}}^{(\theta_0)}(x)]$ is obtained by assembling the element interpolation vectors for $\theta_0$.\\

\noindent\textbf{Expressions for the force vector and mass matrix of the Timoshenko beam f-FEM:}
By using the interpolations for the displacement fields, the virtual work is approximated as:
\begin{equation}
\label{eq: FE_expression_virtual_work}
    \delta \mathcal{V} = \delta \{U\}^T \int_{\Omega_T} \left[ \big\{\mathbb{L}^{(u_0)}\big\}^T F_x + \big\{\mathbb{L}^{(w_0)}\big\}^T F_z + \big\{\mathbb{L}^{(\theta_0)}\big\}^T M_{\theta_0} \right] \mathrm{d}\Omega_T = \delta \{U\}^T [F_T]
\end{equation}
where the row vectors $\big\{\mathbb{L}^{(\square)}\big\} (\square \in \{u_0,w_0,\theta_0\})$ are obtained by assembling the element interpolation vectors given in Eq.~(\ref{eq: u_interpolation}). Similarly, the approximation for the kinetic energy is obtained as:
\begin{equation}
\label{eq: FE_expression_virtual_kinetic_energy}
    \delta T = -\delta \{U\}^T \left[ \int_{\Omega_T} \left\{\bar{\mathbb{L}} \right\} \{I_0, ~I_0, ~I_2\}^T \left\{\bar{\mathbb{L}} \right\}^T \mathrm{d}\Omega_T \right]\{\ddot{U}\} =  -\delta \{U\}^T [M_T] \{\ddot{U}\} 
\end{equation}
where $\left\{\bar{\mathbb{L}} \right\} = \left\{ \big\{\mathbb{L}^{(u_0)}\big\}^T,  \big\{\mathbb{L}^{(w_0)}\big\}^T, \big\{\mathbb{L}^{(\theta_0)}\big\}^T \right\}$.

\section{Comparison of the energy contributions by axial and transverse gradients of axial strain}
\begin{figure}[H]
    \centering
    \includegraphics[width=0.5\textwidth]{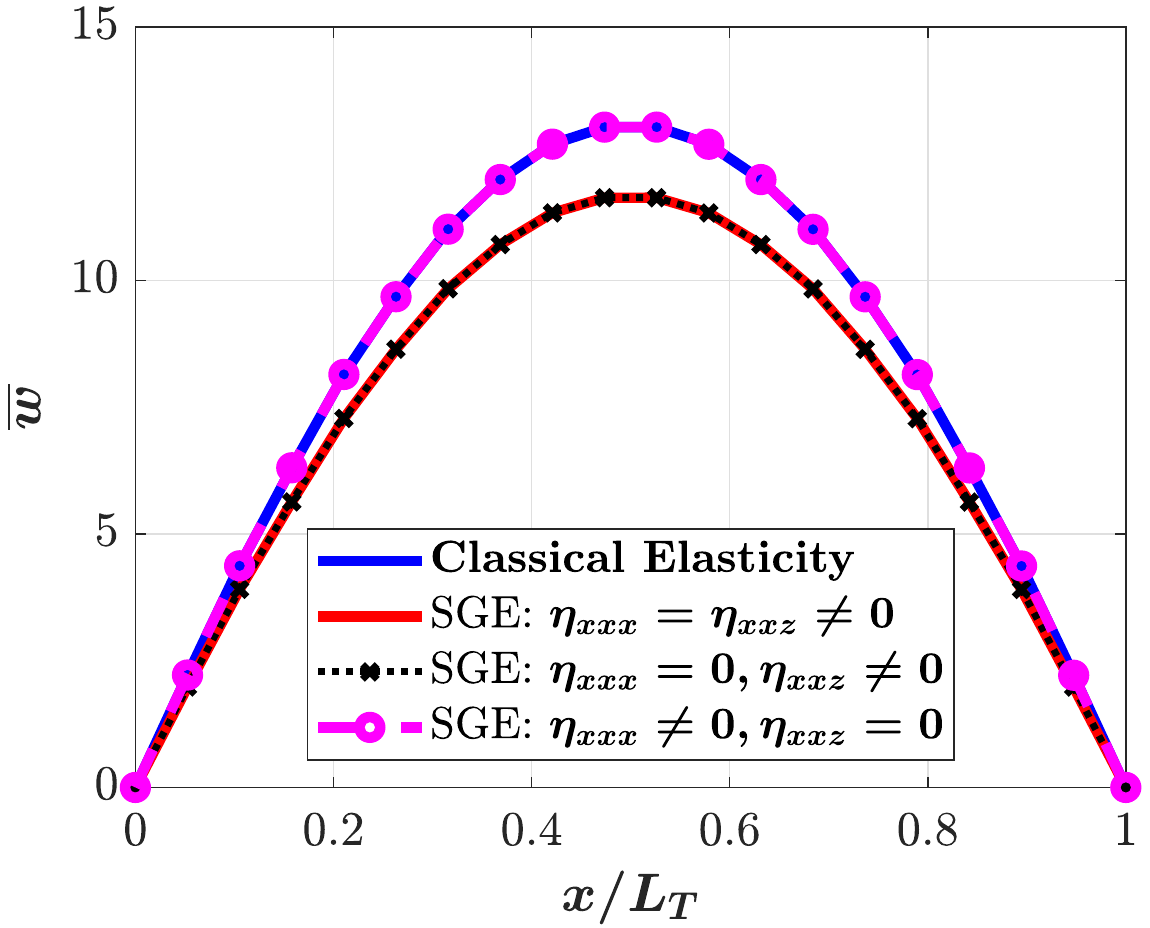}
    \caption{Comparison of the contributions of the axial gradient of the normal strain ($\eta_{xxx}$) and the transverse gradient of the normal strain ($\eta_{xxz}$) to the static response of a beam obtained using Galerkin solutions \cite{sidhardh2018element}.}
    \label{fig: comparison}
\end{figure}

Following the classical first-order strain gradient elasticity, the ratio of the contribution of the axial gradient of the normal strain $(\eta_{xxx})$ and the transverse gradient of the normal strain $(\eta_{xxz})$, to the potential energy of the beam, is obtained as:
\begin{equation}
    \label{eq: beam_CSG_ratio}
    \mathcal{R} = \frac{ E I_{zz} l_*^2 \left( D^2_x \theta_0 \right)^2}{ EA_T l_*^2 \left( D_x^1 \theta_0 \right)^2}
\end{equation}
where $A_T$ and $I_{zz}$ denote the area of cross-section and area moment of inertia of the beam, respectively. Consider the following transformation of the axial variable:
\begin{equation}
    \label{eq: chnage_of_variables}
    x \rightarrow \overline{x}L_T
\end{equation}
Under the above transformation the ratio $\mathcal{R}$ is obtained as:
%\begin{equation}
%    \label{eq: beam_CSG_ratio_modified}
%    \mathcal{R} = \left[ I_{zz} \left( \frac{D^2_{\overline{x}} \theta_0}{L_T^2} \right)^2 \right] \bigg{/} \left[ A_T \left( \frac{D_{\overline{x}}^1 \theta_0}{L_T} \right)^2 \right]
%\end{equation}
%which can be further simplified as:
\begin{equation}
    \label{eq: beam_CSG_ratio_final}
    \mathcal{R} = \frac{1}{12} \left( \frac{h}{L_T} \right)^2 \left[ \frac{D_{\overline{x}}^2 \theta_0}{D_{\overline{x}}^1 \theta_0} \right]^2 \approx O\left( \left(\frac{h}{L_T} \right)^2 \right)
\end{equation}
This indicates that the contribution of the strain-gradient $\eta_{xxx}$ to the potential energy can be ignored compared to $\eta_{xxz}$ for slender beams. This claim is also verified by using the Galerkin solutions provided in \cite{sidhardh2018element} for the static response of beams via first-order strain gradient elasticity. The results are provided in Fig.~(\ref{fig: comparison}). As evident from Fig.~(\ref{fig: comparison}), the static response obtained by ignoring the contribution of $\eta_{xxx}$ to the strain energy, in comparison to $\eta_{xxz}$, closely matches the response obtained without ignoring the same. Note also that, upon ignoring the contribution of $\eta_{xxz}$ to the strain energy, the obtained static response closely matches the classical local elasticity solution. This further indicates that the contribution of $\eta_{xxz}$ to the potential energy is significant when compared to the contribution of $\eta_{xxx}$. 
For the fractional-order formulation, the ratio $\mathcal{R}$ is:
\begin{equation}
    \label{eq: fractional_beam_CSG_ratio_final}
    \mathcal{R} = \frac{1}{12} \left( \frac{h}{L_T} \right)^2 \left[ \frac{D_{\overline{x}}^{\alpha_2} \left( D_{\overline{x}}^{\alpha_1} \theta_0 \right)}{D_{\overline{x}}^{\alpha_1} \theta_0} \right]^2 \approx O\left( \left(\frac{h}{L_T} \right)^2 \right)
\end{equation}
It immediately follows that the previous arguments on the relative energy contributions of $\eta_{xxx}$ and $\eta_{xxz}$ directly extend to the fractional-order formulation.

%\putbib[report]
%\end{bibunit}

\end{document}